\documentclass[11pt,english]{amsart}
\usepackage{amsfonts, amssymb, amsmath, amsthm, eucal, latexsym}
\usepackage{fullpage}

\usepackage[english]{babel}
\usepackage{xy, xypic,amsfonts}
\usepackage{pst-all}
\usepackage{graphicx}
\usepackage{psfrag}  

\theoremstyle{plain}
\newtheorem{theorem}{Theorem}[section]
\newtheorem{lemma}[theorem]{Lemma}
\newtheorem{corollary}[theorem]{Corollary}
\newtheorem{proposition}[theorem]{Proposition}

\theoremstyle{definition}
\newtheorem{definition}[theorem]{Definition}
\newtheorem{example}[theorem]{Example}
\newtheorem{examples}[theorem]{Examples}
\newtheorem{notation}[theorem]{Notation}

\theoremstyle{remark}
\newtheorem{remark}[theorem]{Remark}

\theoremstyle{theorem}

\def\g{{\mathfrak{g}}}  
\def\q{{\mathfrak{q}}}  
\def\l{{\mathfrak{l}}}  
\def\z{{\mathfrak{z}}}  
\def\u{{\mathfrak{u}}}  
\def\a{{\mathfrak{a}}}  %
\def\b{{\mathfrak{b}}}  
\def\t{{\mathfrak{t}}}  
\def\h{{\mathfrak{h}}}  
\def\p{{\mathfrak{p}}}  
\def\m{{\mathfrak{m}}}  
\def\n{{\mathfrak{n}}}  
\def\s{{\mathfrak{s}}}  

\def\C{{\mathbb{C}}}  

\def\uq#1{{u}_{#1}^{-}}

\def\rk{{\rm rk}\hskip .1em} 
\def\ind{\operatorname{ind \, \! }\nolimits} 
\def\kg#1{\mathrm{k}_{#1}} 

\def\varp#1#2{(\varphi_{#1})|_{#2}} 
\def\var#1{\varphi_{#1}} 
\def\eps#1{\varepsilon_{#1}}  
\def\cali#1{\mathcal{#1}}

\title[Quasi-reductive (bi)parabolic subalgebras]{Quasi-reductive (bi)parabolic subalgebras in reductive Lie algebras.}
\author[K. Baur]{Karin Baur}
\address{Karin Baur, ETH Z\"urich\\
Departement Mathematik \\
R\"amistrasse 101\\
8092 Z\"urich, Switzerland}
\email{baur@math.ethz.ch}

\author[A. Moreau]{Anne Moreau}
\address{Anne Moreau, LMA\\
Boulevard Marie et Pierre Curie\\
86962 Futuroscope Chasseneuil Cedex, France}
\email{anne.moreau@math.univ-poitiers.fr}



\begin{document}

\large

\maketitle

\begin{abstract}
We say that a finite dimensional Lie algebra is quasi-reductive if it has a linear form 
whose stabilizer for the coadjoint representation, modulo the center, is a reductive Lie algebra with a center consisting of semisimple elements. 
Parabolic subalgebras of a semisimple Lie algebra are not always quasi-reductive (except in types A or C by work of Panyushev). 
The classification of quasi-reductive parabolic subalgebras in the classical case has been recently achieved in unpublished work of Duflo, 
Khalgui and Torasso. 
In this paper, we investigate the quasi-reductivity of biparabolic subalgebras of reductive Lie algebras. 
Biparabolic (or seaweed) subalgebras are the intersection of two parabolic subalgebras whose sum is the total Lie algebra.
As a main result, we complete the classification of quasi-reductive parabolic subalgebras of reductive Lie algebras by considering the exceptional cases. 
\end{abstract}

{\footnotesize
\begin{quote}
   
   {\sc R\'esum\'e.} \emph{Une alg\`ebre de Lie de dimension finie est dite quasi-r\'eductive si elle poss\`ede une forme lin\'eaire 
   dont le stablisateur pour la repr\'esentation coadjointe, modulo le centre, est une alg\`ebre de Lie r\'eductive 
   avec un centre form\'e d'\'el\'ements semi-simples. 
   Les sous-alg\`ebres paraboliques d'une alg\`ebre de Lie semi-simple ne sont pas toujours quasi-r\'eductives  
   (sauf en types A ou C d'apr\`es un r\'esultat de Panyushev). 
   R\'ecemment, Duflo, Khalgui and Torasso ont termin\'e la classification des sous-alg\`ebres paraboliques quasi-r\'eductives dans le cas classique.
   Dans cet article nous \'etudions la quasi-r\'eductivit\'e des sous-alg\`ebres biparaboliques des alg\`ebres de Lie r\'eductives.
   Les sous-alg\`ebres biparaboliques sont les intersections de deux sous-alg\`ebres paraboliques dont la somme est l'alg\`ebre de Lie ambiante.
   Notre principal r\'esultat est la compl\'etion de la classification
   des sous-alg\`ebres paraboliques quasi-r\'eductives des alg\`ebres de Lie r\'eductives.}
\end{quote}

%
\section*{Introduction}\label{Intro}
%
Let $G$ be a complex connected linear algebraic Lie group.
Denote by $\g$ its Lie algebra.
The group $G$ acts on the dual $\g^*$ of $\g$ by the coadjoint action.
For $f \in \g^*$, we denote by $G(f)$ its stabilizer in $G$; it always contains the center $Z$ of $G$.
One says that a linear form $f\in\g^*$ has \emph{reductive type} if the quotient $G(f)/Z$ is a reductive subgroup of GL$(\g)$.
The Lie algebra $\g$ is called \emph{quasi-reductive} if it has linear forms of reductive type.
This notion goes back to M.~Duflo.
He initiated the study of such Lie algebras because of applications in harmonic
analysis, see \cite{Du1}.
For more details about linear forms of reductive type and quasi-reductive
Lie algebras we refer the reader to Section~\ref{S-1}.

Reductive Lie algebras are obviously quasi-reductive Lie algebras since in that case, $0$
is a linear form of reductive type.
Biparabolic subalgebras form a very interesting class of non-reductive Lie algebras.
They naturally extend the classes of parabolic subalgebras and of Levi subalgebras.
The latter are clearly quasi-reductive since they are reductive subalgebras..
Biparabolic subalgebras were introduced by V.~Dergachev and A.~Kirillov in the case
$\g=\mathfrak{sl}_n$,~see \cite{DK}.
A {\em biparabolic} subalgebra or {\em seaweed} subalgebra (of a semisimple Lie
algebra) is the
intersection of two parabolic subalgebras whose sum is the total Lie algebra.

In this article, we are interested in the classification of quasi-reductive
(bi)parabolic subalgebras.
Note that it is enough to consider the case of (bi)parabolic subalgebras of the simple
Lie algebras,
cf.~Remark~\ref{r-sim}.

In the classical cases, various results are already known:
All biparabolic subalgebras of
$\mathfrak{sl}_n$ and $\mathfrak{sp}_{2n}$ are quasi-reductive as has been proven by D.~Panyushev in~\cite{Pa3}.
The case of orthogonal Lie algebras is more complicated:
On one hand, there are parabolic subalgebras of orthogonal Lie algebras
which are not quasi-reductive, as P.~Tauvel and R.W.T.~Yu have shown (Section 3.2 of~\cite{ty0}).
On the other hand, D.~Panyushev and A.~Dvorsky exhibit many quasi-reductive parabolic subalgebras in~\cite{Dv} and~\cite{Pa3}
by constructing linear forms with the desired properties.
Recently, M.~Duflo, M.S.~Khalgui and P.~Torasso have obtained the classification of quasi-reductive parabolic subalgebras of
the orthogonal Lie algebras in unpublished work,~\cite{DKT}.
They were able to characterize quasi-reductive parabolic subalgebras in terms of the flags stabilized by the subalgebras.\\

The main result of this paper is the completion of the classification of quasi-reductive parabolic subalgebras of simple Lie algebras.
This is done in Section~\ref{S-5} (Theorem~\ref{t-EF} and Theorem~\ref{t-E6}).
Our goal is ultimately to describe all quasi-reductive biparabolic subalgebras.
Thus, in the first sections we present results concerning biparabolic subalgebras to remain in a general setting as far as possible.
For the remainder of the introduction, $\mathfrak{g}$ is a finite dimensional complex semisimple Lie algebra.\\

The paper is organized as follows:\\

In Section~\ref{S-1} we introduce the main notations and definitions.
We also include in this section a short review of known results about biparabolic subalgebras, including the description of
quasi-reductive parabolic subalgebras in the classical Lie algebras (Subsection~\ref{Ss-cla}).
In Section~\ref{S-2}, we describe two methods of reduction, namely the transitivity property
(Theorem~\ref{t-tran}) and the additivity property (Theorem~\ref{t-add}).
As a first step of our classification, we exhibit in Section~\ref{S-3} a large collection of quasi-reductive biparabolic subalgebras of $\g$ (Theorem~\ref{t-lac}).
Next, in Section~\ref{S-4}, we consider \emph{the} non quasi-reductive parabolic subalgebras of $\g$, for simple $\g$ of exceptional
type (Theorems~\ref{t-one},~\ref{t-nE78} and~\ref{t-nE6}).
This is a crucial part of the paper.
Indeed, to study the quasi-reductivity, we can make explicit computations (cf.~Section~\ref{S-5}) while it is much trickier to prove that a Lie algebra is not quasi-reductive.
Using the results of Sections~\ref{S-2},~\ref{S-3} and~\ref{S-4}, we are able to cover
a large number of parabolic subalgebras.
The remaining cases are dealt with in Section~\ref{S-5} (Theorem~\ref{t-2},
Propositions~\ref{p-min} and~\ref{p-GAP}).
This completes the classification of quasi-reductive parabolic subalgebras of $\g$
(Theorems~\ref{t-EF} and~\ref{t-E6}, see also Tables~\ref{T-E8} and~\ref{T-E6}).

At this place, we also want to point out that in~\cite{MY},
O.~Yakimova and the second author study the {\it maximal reductive stabilizers} of
quasi-reductive parabolic subalgebras of $\g$. This piece of work yields an alternative proof
of Proposition~\ref{p-GAP} which is not based on the computer programme \texttt{GAP},
see Remark~\ref{r-yak}.

\setcounter{tocdepth}{1}
\tableofcontents

\noindent
\textbf{Acknowledgment:}
We thank M.~Duflo for introducing us to the subject of quasi-reductive subalgebras.
We also thank P.~Tauvel and R.W.T.~Yu for useful discussions.
Furthermore, we thank W.~de~Graaf and J.~Draisma for helpful hints in the use
of \texttt{GAP}.
At this point, we also want to thank the referee for the very useful comments and
suggestions.

%
\section{Notations, definitions and basic facts} \label{S-1}
%

In this section, we recall a number of known results that will be used in the sequel.

\subsection{}
Let $\g$ be a complex Lie algebra of a connected linear algebraic Lie
group $G$.
Denoting by $\g(f)$ the Lie algebra of $G(f)$,
we have $\g(f)=\{x \in \g \ | \ (\mathrm{ad}^*x)(f)=0 \}$
where $\mathrm{ad}^*$ is the coadjoint representation of $\g$.
Recall that a linear form $f \in \g^*$ is of \emph{reductive type} if $G(f)/Z$ is a reductive Lie subgroup of GL$(\g)$.
We can reformulate this definition as follows:

\begin{definition}
An element $f$ of $\g^*$ is said to be {\it of reductive type}
if $\g(f)/\z$ is a reductive Lie algebra whose center consists of semisimple elements of $\g$
where $\z$ is the center of $\g$.
\end{definition}

Recall that a linear form $f \in \g^*$ is {\em regular} if the dimension of $\g(f)$ is as small as possible.
By definition, the {\em index} of $\g$, denoted by $\ind\g$, is the dimension of the stabilizer of a regular linear form.
The index of various 
special classes of subalgebras of reductive Lie algebras has been studied by
several authors, cf.~\cite{Pa2},~\cite{Ya},~\cite{Mo1},~\cite{Mo2}.
For the index of seaweed algebras, we refer to~\cite{Pa1},~\cite{Dv},~\cite{ty0},~\cite{TY2},~\cite{Jo1} and~\cite{Jo2}.

Recall that $\g$ is called {\it quasi-reductive} if it has linear forms of reductive type.
From Duflo's work~\cite[\S\S I.26-27]{Du1} one deduces the following result
about regular linear forms of reductive type:

\begin{proposition}\label{p-Zar}
Suppose that $\g$ is quasi-reductive.
The set of regular linear forms of reductive type forms a Zariski open dense subset of $\g^*$.
\end{proposition}

\subsection{}\label{Ss-1.2}

From now on, $\mathfrak{g}$ is a complex finite dimensional semisimple
Lie algebra.
The dual of $\g$ is identified with $\g$ through the Killing form of $\g$.
For $u \in \g$, we denote by $\varphi_u$ the corresponding element of $\g^*$.
For $u\in\g$, the restriction of $\varphi_u$, to any subalgebra $\a$ of $\mathfrak{g}$
will be denoted by $\varp{u}{\mathfrak{a}}$.

Denote by $\pi$ the set of simple roots with respect to a fixed triangular decomposition
$$\mathfrak{g}=\mathfrak{n}^+ \oplus \mathfrak{h} \oplus \mathfrak{n}^-$$
of $\mathfrak{g}$, and by $\Delta_{\pi}$ (respectively $\Delta_{\pi}^+$, $\Delta_{\pi}^{-}$)
the corresponding root system (respectively positive root system, negative root system).
If $\pi'$ is a subset of $\pi$, we denote by $\Delta_{\pi'}$ the root subsystem of $\Delta_{\pi}$ generated
by $\pi'$ and we set $\Delta_{\pi'}^{\pm}=\Delta_{\pi'} \cap \Delta_{\pi}^{\pm}$.
For $\alpha \in \Delta_{\pi}$, denote by $\mathfrak{g}_{\alpha}$ the $\alpha$-root subspace of $\g$ and let $h_\alpha$ be the unique element of
$[\mathfrak{g}_{\alpha},\mathfrak{g}_{-\alpha}]$ such that $\alpha(h_{\alpha})=2$.
For each $\alpha \in \Delta_{\pi}$, fix $x_{\alpha} \in \g_{\alpha}$ so that the family
$\{x_{\alpha}, h_{\beta} \ ; \  \alpha \in \Delta_{\pi}, \beta \in \pi \}$ is a Chevalley basis of $\g$.
In particular, for non-colinear roots $\alpha$ and $\beta$,
we have $[x_{\alpha},x_{\beta}] =\pm(p+1)x_{\alpha+\beta}$ if $\beta-p\alpha$ is the source of the $\alpha$-string through $\beta$.

We briefly recall a classical construction due to B.~Kostant.
It associates to a subset of $\pi$ a system of strongly orthogonal positive roots in $\Delta_{\pi}$.
This construction is known to be very helpful to obtain regular forms on biparabolic subalgebras of $\mathfrak{g}$.
For a recent account about the \emph{cascade} construction of Kostant, we refer to~\cite[\S1.5]{TY2} or~\cite[\S40.5]{TY1}.

For $\lambda$ in $\mathfrak{h}^{*}$ and $\alpha \in \Delta_{\pi}$, we shall write $\langle
\lambda, \alpha^{\vee} \rangle$ for $\lambda(h_{\alpha})$.
Recall that two roots $\alpha$ and $\beta$ in $\Delta_{\pi}$ are said to be {\it strongly orthogonal} if neither $\alpha+\beta$ nor $\alpha -\beta$ is in  $\Delta_{\pi}$.
Let $\pi'$ be a subset of $\pi$.
The {\em cascade} $\mathcal{K}_{\pi'}$ of $\pi'$ is defined by induction on the cardinality of $\pi'$ as follows:
\begin{itemize}
\item[(1)] $\mathcal{K}(\emptyset)=\emptyset$,
\item[(2)] If $\pi'_1$,\ldots,$\pi'_r$ are the connected components of $\pi'$, then~
$\mathcal{K}_{\pi'}=\mathcal{K}_{\pi'_1} \cup \cdots \cup \mathcal{K}_{\pi'_r}$,
\item[(3)] If $\pi'$ is connected,  then $\mathcal{K}_{\pi'}=\{\pi'\}
\cup \mathcal{K}_{T}$ where $T=\{ \alpha \in \pi' \ | \ \langle \alpha, \eps{\pi'}^{\vee} \rangle =0 \}$
and $\eps{\pi'}$ is the highest positive root of $\Delta_{\pi'}^+$.
\end{itemize}
For $K \in \mathcal{K}_{\pi'}$, set
$$\Gamma_{K}=\{ \alpha \in \Delta_{K} \ | \
 \langle \alpha, \varepsilon_{K}^{\vee} \rangle > 0 \} \quad \textrm{and} \quad
\Gamma^{0}_{K}=
\Gamma_{K} \setminus\{\varepsilon_{K}\}\, .$$
Notice that the subspace $\sum\limits_{K\in\Gamma_K} \g_{\alpha}$ is a Heisenberg Lie algebra
whose center is $\g_{\varepsilon_K}$.

The cardinality $\kg{\pi}$ of $\mathcal{K}_{\pi}$ only depends on $\mathfrak{g}$;
it is independent of the choices of $\mathfrak{h}$ and $\pi$.
The values of $\kg{\pi}$
for the different types of simple Lie algebras are given in Table~\ref{Tkg};
in this table, for a real number $x$, we denote by $[x]$ the largest integer
$\leq x$.

{
\begin{table}[h] \tiny
\begin{center}
\begin{tabular}{|c|c|c|c|c|c|c|c|c|}
\hline
  & & & & & & & & \\
$\mathrm{A}_{\ell}, \ell \geq 1$  &
$\mathrm{B}_{\ell}, \ell \geq 2$  &
$\mathrm{C}_{\ell}, \ell \geq 3$  &
$\mathrm{D}_{\ell}, \ell \geq 4$  &
$\mathrm{G}_2$  &
$\mathrm{F}_4$  &
$\mathrm{E}_6$  &
$\mathrm{E}_7$  &
$\mathrm{E}_8$   \\
  & & & & & & & & \\
\hline
 & & & & & & & & \\
$\left[ \displaystyle{\frac{\ell+1}{2}} \right]$ & $\ell$ & $\ell$ &
$2 \left[ \displaystyle{\frac{\ell}{2}} \right]$ &
$2$ & $4$ & $4$ & $7$ & $8$ \\
  & & & & & & & & \\
\hline
\end{tabular}
\vspace{.5cm}
\caption{\label{Tkg} $\kg{\pi}$ for the simple  Lie algebras.}
\end{center}
\end{table}}

For $\pi'$ a subset of $\pi$, we denote by $\cali{E}_{\pi'}$ the set of
the highest roots $\eps{K}$ where $K$ runs over the elements of the cascade of $\pi'$.
By construction, the subset $\cali{E}_{\pi'}$ is a family of pairwise strongly orthogonal roots in
$\Delta_{\pi'}$.
For the convenience of the reader, the set $\cali{E}_{\pi}$, for each simple Lie algebra of type $\pi$,
is described in the Tables \ref{T-cla} and \ref{T-exc}.
We denote by $E_{\pi'}$ the subspace of $\mathfrak{h}^*$ which is generated
by the elements of $\cali{E}_{\pi'}$.

{\begin{table}[h]\tiny
\begin{tabular}{lcl}
\hline
$\mathrm{A}_{\ell}$, \ $\ell \geq 1$: &
\begin{pspicture}(-1.5,-0)(3.5,.75)
\pscircle(-1,0){1mm}
\pscircle(0,0){1mm}
\pscircle(1,0){1mm}
\pscircle(2,0){1mm}
\pscircle(3,0){1mm}
\psline(-0.1,0)(-0.9,0)
\psline[linestyle=dotted](0.1,0)(0.9,0)
\psline(1.1,0)(1.9,0)
\psline(2.1,0)(2.9,0)
\rput[b](-1,0.2){$\alpha_1$}
\rput[b](0,0.2){$\alpha_2$}
\rput[b](2,0.2){$\alpha_{\ell-1}$}
\rput[b](3,0.2){$\alpha_{\ell}$}
\end{pspicture}
& $\{ \eps{i}=\alpha_i + \cdots + \alpha_{i+(\ell-2i +1)}, \ i \leq
\left[ \displaystyle{\frac{\ell+1}{2}}\right] \}$\\
$\mathrm{B}_{\ell}$, \ $\ell \geq 2$ :  &
\begin{pspicture}(-1.5,-0)(3.5,.75) 
\pscircle(-1,0){1mm}
\pscircle(0,0){1mm}
\pscircle(1,0){1mm}
\pscircle(2,0){1mm}
\pscircle(3,0){1mm}
\psline(-0.1,0)(-0.9,0)
\psline[linestyle=dotted](0.1,0)(0.9,0)
\psline(1.1,0)(1.9,0)
\psline(2.09,0.05)(2.91,0.05)
\psline(2.09,-0.05)(2.91,-0.05)
\rput[b](-1,0.2){$\alpha_1$} 
\rput[b](0,0.2){$\alpha_2$}  
\rput[b](2,0.2){$\alpha_{\ell-1}$}
\rput[b](3,0.2){$\alpha_{\ell}$}
\rput(2.5,0){$>$}
\end{pspicture}
&$\{ \eps{i}=\alpha_{i-1} + 2 \alpha_i + \cdots + 2 \alpha_{\ell} , \ i \textrm{ even}, \
i \leq \ell \} \cup  \{ \eps{i}=\alpha_i  ,  \ i \textrm{ odd},  i \leq \ell \}$\\
&& \\
$\mathrm{C}_{\ell}$, $\ell \geq 3$: &
\begin{pspicture}(-1.5,-0)(3.5,.75)
\pscircle(-1,0){1mm}
\pscircle(0,0){1mm}
\pscircle(1,0){1mm}
\pscircle(2,0){1mm}
\pscircle(3,0){1mm}
\psline(-0.1,0)(-0.9,0)
\psline[linestyle=dotted](0.1,0)(0.9,0)
\psline(1.1,0)(1.9,0)
\psline(2.09,0.05)(2.91,0.05)
\psline(2.09,-0.05)(2.91,-0.05)
\rput[b](-1,0.2){$\alpha_1$} 
\rput[b](0,0.2){$\alpha_2$} 
\rput[b](2,0.2){$\alpha_{\ell-1}$}
\rput[b](3,0.2){$\alpha_\ell$}
\rput(2.5,0){$<$}
\end{pspicture}
& $\{ \eps{i}=2 \alpha_i + \cdots +2 \alpha_{\ell-1} + \alpha_{\ell} , \
 i \leq \ell-1 \} \cup  \{ \eps{\ell} = \alpha_{\ell}  \}$\\
 & & \\
$\mathrm{D}_{\ell}$, $\ell$ even, $\ell \geq 4$: &
\begin{pspicture}(-1.5,-.5)(3.6,0.3) 
\pscircle(-1,-1){1mm}
\pscircle(0,-1){1mm}
\pscircle(1,-1){1mm}
\pscircle(2,-1){1mm}
\pscircle(3,-0.3){1mm}
\pscircle(3,-1.7){1mm}
\psline(-0.1,-1)(-0.9,-1)
\psline[linestyle=dotted](0.1,-1)(0.9,-1)
\psline(1.1,-1)(1.9,-1)
\psline(2.07,-.95)(2.92,-.32)
\psline(2.07,-1.05)(2.92,-1.68)
\rput[b](-1,-0.8){$\alpha_1$} 
\rput[b](0,-0.8){$\alpha_2$} 
\rput[b](2,-0.8){$\alpha_{\ell-2}$}
\rput[l](3.1,-0.2){$\alpha_{\ell}$}
\rput[l](3.1,-1.8){$\alpha_{\ell-1}$}
\end{pspicture}
&   $\{ \eps{i}=\alpha_{i-1} + 2 \alpha_{i} + \cdots +2 \alpha_{\ell-2} +
   \alpha_{\ell-1} +\alpha_{\ell}, \ i \textrm{ even}, \  i < \ell-1 \}$ \\
& & \\
 & & $\quad\quad \cup \  \{ \eps{i}=\alpha_i , \ i \textrm{ odd}, \  i < \ell \} \  \cup \
 \{ \eps{\ell}=\alpha_{\ell} \}$ \\
  &&\\

\medskip
$\mathrm{D}_{\ell}$, $\ell$ odd, $\ell \geq 5$: &
\begin{pspicture}(-1.5,-0.5)(3.6,.3)
\pscircle(-1,-1){1mm}
\pscircle(0,-1){1mm}
\pscircle(1,-1){1mm}
\pscircle(2,-1){1mm}
\pscircle(3,-0.3){1mm}
\pscircle(3,-1.7){1mm}
\psline(-0.1,-1)(-0.9,-1)
\psline[linestyle=dotted](0.1,-1)(0.9,-1)
\psline(1.1,-1)(1.9,-1)
\psline(2.07,-.95)(2.92,-.32)
\psline(2.07,-1.05)(2.92,-1.68)
\rput[b](-1,-0.8){$\alpha_1$} 
\rput[b](0,-0.8){$\alpha_2$} 
\rput[b](2,-0.8){$\alpha_{\ell-2}$}
\rput[l](3.1,-0.2){$\alpha_{\ell}$}
\rput[l](3.1,-1.8){$\alpha_{\ell-1}$}
\end{pspicture}
 & $\{ \eps{i}=\alpha_{i-1} + 2 \alpha_{i} + \cdots +2 \alpha_{\ell-2} +
  \alpha_{\ell-1} +\alpha_{\ell} , \ i \textrm{ even }, \ i < \ell-1  \}$ \\
 & & \\
& & $\quad\quad \cup \ \{\eps{i}=\alpha_i , \ i \textrm{ odd }, \  i < \ell \} \ \cup  \ \{
  \eps{\ell-1}=\alpha_{\ell-2} + \alpha_{\ell-1} +\alpha_{\ell}\}$\\
 &&\\
 &&\\
\hline
\end{tabular}
\vspace{.01cm}
\caption{\label{T-cla} $\cali{E}_{\pi}$ for the classical Lie algebras.}
\end{table}}

{\begin{table}[h]\tiny

\begin{tabular}{cll}
\hline
$\mathrm{G}_2$: &
\begin{pspicture}(-1.2,-0.05)(2,0.4)
\pscircle(-1,0){1mm}
\pscircle(0,0){1mm}
\psline(-0.09,0.06)(-0.91,0.06)
\psline(-0.09,0)(-0.91,0)
\psline(-0.09,-0.06)(-0.91,-0.06)
\rput[b](-1,0.15){$\alpha_{1}$}
\rput[b](0,0.15){$\alpha_{2}$}
\rput(-0.5,0){$>$}
\end{pspicture}
& $\{\eps{1}=\begin{array}{c}23\end{array}, \ \eps{2}=\begin{array}{c}01\\\end{array}\}$\\

\medskip
$\mathrm{F}_4$: &
\begin{pspicture}(-2.2,-0.05)(2,0.4)
\pscircle(-2,0){1mm}
\pscircle(-1,0){1mm}
\pscircle(0,0){1mm}
\pscircle(1,0){1mm}
\psline(-1.9,0)(-1.1,0)
\psline(-0.09,0.05)(-0.91,0.05)
\psline(-0.09,-0.05)(-0.91,-0.05)
\psline(0.1,0)(0.9,0)
\rput[b](-2,0.15){$\alpha_{1}$}
\rput[b](-1,0.15){$\alpha_{2}$}
\rput[b](0,0.15){$\alpha_{3}$}
\rput[b](1,0.15){$\alpha_{4}$}
\rput(-0.5,0){$>$}
\end{pspicture}
&$\{\eps{1}=\begin{array}{c}2342\\\end{array}, \
\eps{2}=\begin{array}{c}0122\end{array}, \
\eps{3}=\begin{array}{c}0120 \end{array}, \
\eps{4}=\begin{array}{c}0100\end{array}\}$\\

\medskip
$\mathrm{E}_6$: &
\begin{pspicture}(-2.2,-0.05)(2,0.3)
\pscircle(-2,0){1mm}
\pscircle(-1,0){1mm}
\pscircle(0,0){1mm}
\pscircle(1,0){1mm}
\pscircle(2,0){1mm}
\pscircle(0,-0.5){1mm}
\psline(-1.9,0)(-1.1,0)
\psline(-0.1,0)(-0.9,0)
\psline(0.1,0)(0.9,0)
\psline(1.1,0)(1.9,0)
\psline(0,-0.1)(0,-0.4)
\rput[b](-2,0.15){$\alpha_1$}
\rput[b](-1,0.15){$\alpha_3$}
\rput[b](0,0.15){$\alpha_4$}
\rput[b](1,0.15){$\alpha_5$}
\rput[b](2,0.15){$\alpha_6$}
\rput[r](-0.2,-0.5){$\alpha_2$}
\end{pspicture}
& $\{\eps{1}=\begin{array}{c}\\12321\\2\end{array}, \
\eps{2}=\begin{array}{c}\\11111\\0\end{array}, \
\eps{3}=\begin{array}{c}\\01110\\0\end{array}, \
\eps{4}=\begin{array}{c}\\00100\\0\end{array}\}$\\

\medskip
$\mathrm{E}_7$: &
\begin{pspicture}(-2.2,-0.05)(3,0.3)
\psline(-1.9,0)(-1.1,0)
\psline(-0.1,0)(-0.9,0)
\psline(0.1,0)(0.9,0)
\psline(1.1,0)(1.9,0)
\psline(2.1,0)(2.9,0)
\psline(0,-0.1)(0,-0.4)
\rput[b](-2,0.15){$\alpha_1$}
\rput[b](-1,0.15){$\alpha_3$}
\rput[b](0,0.15){$\alpha_4$}
\rput[b](1,0.15){$\alpha_5$}
\rput[b](2,0.15){$\alpha_6$}
\rput[b](3,0.15){$\alpha_7$}
\rput[r](-0.2,-0.5){$\alpha_2$}
\rput[b](-2,0.15){$\alpha_1$}
\pscircle(-2,0){1mm}
\pscircle(-1,0){1mm}
\pscircle(0,0){1mm}
\pscircle(1,0){1mm}
\pscircle(2,0){1mm}
\pscircle(3,0){1mm}
\pscircle(0,-0.5){1mm}
\psline(-1.9,0)(-1.1,0)
\psline(-0.1,0)(-0.9,0)
\psline(0.1,0)(0.9,0)
\psline(1.1,0)(1.9,0)
\psline(2.1,0)(2.9,0)
\psline(0,-0.1)(0,-0.4)
\end{pspicture}
& $\{\eps{1}=\begin{array}{c}\\234321\\ \hspace{-0.15cm}2\end{array}, \
\eps{2}=\begin{array}{c}\\012221\\ \hspace{-0.15cm}1\end{array}, \
\eps{3}=\begin{array}{c}\\012100\\ \hspace{-0.15cm}1\end{array}, \
\eps{4}=\begin{array}{c}\\000001\\ \hspace{-0.15cm}0\end{array}$, \\
&&\hspace{0.05cm} $\eps{5}=\begin{array}{c}\\000000\\ \hspace{-0.15cm}1\end{array}, \
\eps{6}=\begin{array}{c}\\010000\\ \hspace{-0.15cm}0\end{array},\
\eps{7}=\begin{array}{c}\\000100\\ \hspace{-0.15cm}0\end{array}\}$\\

\medskip
$\mathrm{E}_8$: &
\begin{pspicture}(-2.05,-0.05)(4.55,0.3)
\pscircle(-2,0){1mm}
\pscircle(-1,0){1mm}
\pscircle(0,0){1mm}
\pscircle(1,0){1mm}
\pscircle(2,0){1mm}
\pscircle(3,0){1mm}
\pscircle(4,0){1mm}
\pscircle(0,-0.5){1mm}
\psline(-1.9,0)(-1.1,0)
\psline(-0.1,0)(-0.9,0)
\psline(0.1,0)(0.9,0)
\psline(1.1,0)(1.9,0)
\psline(2.1,0)(2.9,0)
\psline(3.1,0)(3.9,0)
\psline(0,-0.1)(0,-0.4)
\rput[b](-2,0.15){$\alpha_1$}
\rput[b](-1,0.15){$\alpha_3$}
\rput[b](0,0.15){$\alpha_4$}
\rput[b](1,0.15){$\alpha_5$}
\rput[b](2,0.15){$\alpha_6$}
\rput[b](3,0.15){$\alpha_7$}
\rput[b](4,0.15){$\alpha_8$}
\rput[r](-0.2,-0.5){$\alpha_2$}
\end{pspicture}
& $\{\eps{1}=\begin{array}{c}\\2465432\\ \hspace{-0.25cm}3\end{array},\
\eps{2}=\begin{array}{c}\\2343210\\ \hspace{-0.25cm}2\end{array},\
\eps{3}=\begin{array}{c}\\0122210\\ \hspace{-0.25cm}1\end{array}, \
\eps{4}=\begin{array}{c}\\0121000\\ \hspace{-0.25cm}1\end{array}$,\\
&&\hspace{0.05cm} $\eps{5}=\begin{array}{c}\\0000010\\ \hspace{-0.25cm}0\end{array}, \
\eps{6}=\begin{array}{c}\\0000000\\ \hspace{-0.25cm}1\end{array}, \
\eps{7}=\begin{array}{c}\\0100000\\ \hspace{-0.25cm}0\end{array}, \
\eps{8}=\begin{array}{c}\\0001000\\ \hspace{-0.25cm}0\end{array}\}$\\
\hline
\end{tabular}
\vspace{.01cm}
\caption{\label{T-exc} $\cali{E}_{\pi}$ for the exceptional Lie algebras.}
\end{table}}

%
\subsection{} \label{Ss-bip}
%

A {\it biparabolic} subalgebra of $\mathfrak{g}$ is defined to be the intersection of two parabolic subalgebras whose sum is $\mathfrak{g}$.
This class of algebras has first been studied in the case
of $\mathfrak{sl}_n$ by Dergachev and Kirillov~\cite{DK} under the name of \emph{seaweed} algebras.

For a subset $\pi'$ of $\pi$, we denote by
$\mathfrak{p}_{\pi'}^{+}$ the standard parabolic subalgebra of $\mathfrak{g}$ which  is the subalgebra  generated by $\mathfrak{b}^+=\mathfrak{h}
\oplus \mathfrak{n}^+$ and by $\mathfrak{g}_{-\alpha}$, for $\alpha \in \pi'$.
We denote by  $\mathfrak{p}_{\pi'}^{-}$ the ``opposite parabolic subalgebra'' generated by $\mathfrak{b}^- = \mathfrak{n}^- \oplus \mathfrak{h}$ and
by $\mathfrak{g}_{\alpha}$, for $\alpha \in \pi'$.
Set
$\mathfrak{l}_{\pi'}=\mathfrak{p}_{\pi'}^{+} \cap \mathfrak{p}_{\pi'}^{-}$.
Then $\mathfrak{l}_{\pi'}$ is a  Levi factor of both
parabolic subalgebras $\mathfrak{p}_{\pi'}^{+}$
and $\mathfrak{p}_{\pi'}^{-}$ and we can write
$\mathfrak{l}_{\pi'}=\mathfrak{n}_{\pi'}^{+} \oplus \mathfrak{h} \oplus \mathfrak{n}_{\pi'}^{-}$
where $\mathfrak{n}_{\pi'}^{\pm}= \mathfrak{n}^{\pm} \cap \mathfrak{l}_{\pi'}$.
Let $\mathfrak{m}_{\pi'}^{+}$ (respectively $\mathfrak{m}_{\pi'}^{-}$) be the nilradical of $\mathfrak{p}_{\pi'}^{+}$
(respectively $\mathfrak{p}_{\pi'}^{-}$).
We denote by $\g_{\pi'}$ the derived Lie algebra of $\mathfrak{l}_{\pi'}$ and by
$\z(\l_{\pi'})$ the center of $\l_{\pi'}$.
The Cartan subalgebra $\h \cap \g_{\pi'}$ of $\g_{\pi'}$ will be denoted by $\mathfrak{h}_{\pi'}$.

\begin{definition}
The subalgebra $\mathfrak{q}_{\pi_1,\pi_2}$ of $\mathfrak{g}$
given as follows by the subsets
$\pi_1,\pi_2 \subset \pi$
$$\mathfrak{q}_{\pi_1,\pi_2}
 : = \mathfrak{p}_{\pi_1}^{+} \cap \mathfrak{p}_{\pi_2}^{-}=\mathfrak{n}_{\pi_2}^{+}
  \oplus \mathfrak{h} \oplus  \mathfrak{n}_{\pi_1}^{-}$$
is called the {\em standard} biparabolic subalgebra
(associated to $\pi_1$ and $\pi_2$).
Its nilpotent radical is
$\mathfrak{u}_{\pi_1,\pi_2}:=(\mathfrak{n}_{\pi_2}^{+} \cap \mathfrak{m}_{\pi_1}^{+} )
\oplus (\mathfrak{n}_{\pi_1}^{-} \cap \mathfrak{m}_{\pi_2}^{-})$ and
$\mathfrak{l}_{\pi_1,\pi_2}
  := \mathfrak{l}_{\pi_1 \cap \pi_2}$
is the standard Levi factor of $\mathfrak{q}_{\pi_1,\pi_2}$.
\end{definition}

Any biparabolic subalgebra is conjugate to a standard one, see~\cite[\S 2.3]{TY2} or~\cite[\S 2.5]{Jo1}.
So, for our purpose, it will be enough to consider standard biparabolic subalgebras.

\begin{remark}\label{r-sim}
The classification of quasi-reductive (bi)parabolic subalgebras of reductive Lie algebras can be deduced from the classification
of quasi-reductive (bi)parabolic subalgebras of simple Lie algebras:
A stabilizer of a linear form on $\mathfrak{g}$ is the product of its components on each
of the simple factors of $\mathfrak{g}$ and of
the center of $\mathfrak{g}$.
As a consequence, we may assume that $\g$ is simple without loss of generality.
\end{remark}

Let $\pi_1,\pi_2$ be two subsets of $\pi$.
The dual of $\mathfrak{q}_{\pi_1,\pi_2}$ is identified to
$\mathfrak{q}_{\pi_2,\pi_1}$ via the Killing form of $\mathfrak{g}$.
For $\underline{a}=(a_K)_{K \in \mathcal{K}_{\pi_2}} \in (\C^*)^{\kg{\pi_2}}$
and $\underline{b}=(b_L)_{L \in \mathcal{K}_{\pi_1}} \in (\C^*)^{\kg{\pi_1}}$, set
$$
u(\underline{a},\underline{b})=\sum\limits_{K \in \mathcal{K}_{\pi_2}} a_K x_{-\eps{K}}
+ \sum\limits_{L \in \mathcal{K}_{\pi_1}} b_L x_{\eps{L}}
$$
It is an element of $\u_{\pi_2,\pi_1}$ and the linear form $\varp{u}{\mathfrak{q}_{\pi_1,\pi_2}}$
is a regular element of $\q_{\pi_1,\pi_2}^*$
for any $(\underline{a},\underline{b})$ running through
a nonempty open subset of $(\C^*)^{\kg{\pi_2}+\kg{\pi_1}}$,
cf.~\cite[Lemma 3.9]{TY2}.

We denote by $E_{\pi_1,\pi_2}$ the subspace generated by the elements $\eps{K}$,
for $K \in \mathcal{K}_{\pi_1} \cup \mathcal{K}_{\pi_2} $.
Thus, $\dim E_{\pi_1,\pi_2} =\kg{\pi_1}+\kg{\pi_2} - \dim(E_{\pi_1} \cap E_{\pi_2})$.
As it has been proved in~\cite[\S 7.16]{Jo1}, we have
\begin{eqnarray} \label{eq-ind}
\ind \q_{\pi_1,\pi_2}=(\rk\g-\dim E_{\pi_1,\pi_2})+(\kg{\pi_1}+\kg{\pi_2}-\dim E_{\pi_1,\pi_2})
\end{eqnarray}

\begin{remark} \label{r-0}
By~(\ref{eq-ind}),
the index of $\q_{\pi_1,\pi_2}$ is zero
if and only if $ E_{\pi_1} \cap E_{\pi_2}=\{0\}$ and
$\kg{\pi_1}+\kg{\pi_2}=\rk \g$.
For example, in type $\mathrm{E}_6$, there are exactly fourteen
standard parabolic subalgebras $\p_{\pi'}^+$ with index zero.
The corresponding subsets $\pi'\subset \pi$ of the simple roots are the following:
$$\begin{array}{l}
\{\alpha_1,\alpha_5\} ; \{\alpha_3,\alpha_6\} ;
\{\alpha_1,\alpha_4,\alpha_5\} ; \{\alpha_3,\alpha_4,\alpha_6\} ;
\{\alpha_1,\alpha_5,\alpha_6\} ;\\
\{\alpha_1,\alpha_3,\alpha_6\} ;
\{\alpha_1,\alpha_3,\alpha_5\} ;  \{\alpha_3,\alpha_5,\alpha_6\} ;
\{\alpha_1,\alpha_3,\alpha_4\} ;  \{\alpha_4,\alpha_5,\alpha_6\} ; \\
\{\alpha_1,\alpha_3,\alpha_4,\alpha_5\} ;  \{\alpha_3,\alpha_4,\alpha_5,\alpha_6\} ;
\{\alpha_1,\alpha_2,\alpha_3,\alpha_4\} ; \{\alpha_2,\alpha_4,\alpha_5,\alpha_6\} .
\end{array}$$
This was already observed in the unpublished work~\cite{El} of A.~Elashvili (with a small error).
\end{remark}

In the sequel, we will often make use of the following element of $\u_{\pi_2,\pi_1}$ on our way to
construct reductive forms:
\begin{eqnarray*} \label{eqU1}
\uq{\pi_1,\pi_2} &=& \sum\limits_{\eps{} \in \cali{E}_{\pi_2}\mbox{, }\eps{}\notin \Delta_{\pi_1}^{+}} x_{-\eps{}}
\end{eqnarray*}
If $\pi_2=\pi$, we simply write
$\uq{\pi_1}$ for $\uq{\pi_1,\pi}$ and,
in the special case of $\pi_1=\emptyset$ and $\pi_2=\pi$, we write
$\uq{}$ for $\uq{\emptyset}$.
Let $B$ be the Borel subgroup of $G$ whose Lie algebra
is $\b^+$. We summarize in the following proposition
useful results of Kostant concerning the linear form $\varp{\uq{}}{\b^+}$.
They can be found
in~\cite[Proposition 40.6.3]{TY1}.

\begin{proposition} \label{p-bor}
{\rm (i)} The linear form $\varp{\uq{}}{\b^+}$ is of reductive type for $\b^+$.
More precisely, the stabilizer of $\var{\uq{}}$ in $\b^+$
is the subspace $\bigcap\limits_{K \in \cali{K}_{\pi}} \ker \eps{K}$ of $\h$
of dimension $\rk\g-\kg{\pi}$.

{\rm (ii)} Let $\mathfrak{m}$ be an ideal of $\b^+$ contained in $\mathfrak{n}^+$.
The $B$-orbit of $\varp{\uq{}}{\m}$ in $\mathfrak{m}^*$ is an open dense
subset of $\mathfrak{m}^*$.
\end{proposition}

%
\subsection{} \label{Ss-cla}
%

We end the section by reviewing what is known in the classical case.
First recall that the biparabolic subalgebras of simple Lie algebras
of type A and C are always quasi-reductive as has been shown by D.~Panyushev in~\cite{Pa3}.

The classification of quasi-reductive parabolic subalgebras of the
orthogonal Lie algebras is given in the recent work~\cite{DKT} of Duflo, Khalgui and Torrasso.
Since we will use this result repeatedly, we state it below.

Let $E$ be a complex vector space of dimension $N$ endowed with a nondegenerate symmetric bilinear form.
Denote by $\mathfrak{so}_N$ the Lie algebra of the corresponding orthogonal group.
Let $\mathcal{V}=\{\{0\}=V_0 \subsetneq V_{1}\subsetneq\dots\subsetneq V_{s}=V\}$ be a
flag of isotropic subspaces in $E$,
with $s \ge 1$.
Its stabilizer in $\mathfrak{so}_N$ is a parabolic subalgebra of $\mathfrak{so}_N$ and
any parabolic subalgebra of $\mathfrak{so}_N$ is obtained in this way.
We denote by $\p_{\mathcal{V}}$ the stabilizer of $\mathcal{V}$ in $\mathfrak{so}_N$.

\begin{theorem}\cite{DKT} \label{t-dkt}
Let $\mathcal{V}=\{\{0\}=V_0 \subsetneq V_{1}\subsetneq\dots\subsetneq V_{s}=V\}$ be a flag of isotropic subspaces in $E$ with $s\ge 1$.
Denote by $\mathcal{V}'$ the flag of isotropic subspaces in  $E$ which is equal to $\mathcal{V} \setminus \{V\}$ if $\dim V$ is odd and equal to $N/2$, and equal to $\mathcal{V}$ otherwise.

The Lie algebra $\p_{\mathcal{V}}$ is quasi-reductive if and only if the sequence $\mathcal{V}'$ does not contain two consecutive subspaces of odd dimension.
\end{theorem}

\begin{example} \label{ex-d6}
For $\mathfrak{g}=$D$_6$ there are twelve standard parabolic subalgebras
$\p=\p_{\pi'}^+$ which are not quasi-reductive.
The corresponding subsets $\pi'\subset\pi$ of the simple
roots are the following:
$$
\begin{array}{l}
\{\alpha_2\}, \{\alpha_4\}, \{\alpha_1,\alpha_4\}, \{\alpha_2,\alpha_4\},
\{\alpha_2,\alpha_5\},
\{\alpha_2,\alpha_6\}, \\
\{\alpha_1,\alpha_2,\alpha_4\},
\{\alpha_2,\alpha_3,\alpha_4\},
\{\alpha_2,\alpha_4,\alpha_5\},
\{\alpha_2,\alpha_4,\alpha_6\}, \\
\{\alpha_2,\alpha_5,\alpha_6\},
\{\alpha_2,\alpha_4,\alpha_5,\alpha_6\} .
\end{array}
$$
Among these, the 
connected $\pi'$ are $\{\alpha_2\},  \{\alpha_4\}, \{\alpha_2,\alpha_3,\alpha_4\}$.
\end{example}

Thus it remains to determine the quasi-reductive parabolic subalgebras of
the exceptional Lie algebras. This is our goal.

\section{Methods of reduction} \label{S-2}

In this section, we
develop methods of reduction to deduce the quasi-reductivity
of a parabolic subalgebra from the quasi-reductivity
of other subalgebras.
We assume that $\pi_2=\pi$.
Nevertheless we keep the notations of biparabolic subalgebras
where it is convenient.

\subsection{}

The following theorem seems to be standard. As there is no proof to our knowledge,
we give a short proof here:

\begin{theorem}[Transitivity] \label{t-tran}
Let $\pi'',\pi'$
be subsets of $\pi$ with $\pi''\subset \pi'$.
Suppose that $\cali{K}_{\pi'} \subset \cali{K}_{\pi}$.
Then, $\q_{\pi'',\pi}$ is quasi-reductive
if and only if
$\q_{\pi'',\pi'}$ is.
\end{theorem}

\begin{proof}
Note that the assumption $\cali{K}_{\pi'} \subset \cali{K}_{\pi}$
implies $\ind \q_{\pi'',\pi'}=\ind \q_{\pi'',\pi}+(\kg{\pi}-\kg{\pi'})$
by formula~(\ref{eq-ind}).
Since $\u_{\pi',\pi}$ is an ideal of $\b^+$ contained in $\n^+$, Proposition~\ref{p-bor}(ii) enables to
choose $w'$ in $\l_{\pi'}$ such that both $\varp{w'+\uq{\pi'}}{\q_{\pi'',\pi}}$
and $\varp{w'}{\q_{\pi'',\pi'}}$ are regular linear forms of $\q_{\pi'',\pi}$ and
$\q_{\pi'',\pi'}$ respectively.
Then one can show that $\q_{\pi'',\pi'}(\var{w'})=\q_{\pi'',\pi}(\var{w'+\uq{\pi'}})
\oplus \sum\limits_{K \in \cali{K}_{\pi}\setminus \cali{K}_{\pi'}} \C h_{\eps{K}}$.
By Proposition~\ref{p-Zar}, if $\q_{\pi'',\pi}$ (respectively $\q_{\pi'',\pi'}$) is quasi-reductive,
then we can assume furthermore that $\varp{w'+\uq{\pi'}}{\q_{\pi'',\pi}}$  (respectively $\varp{w'}{\q_{\pi'',\pi'}}$) has reductive type.
Hence the equivalence of the theorem follows.
\end{proof}

Suppose that $\g$ is simple and
let $\widetilde{\pi}$ be the subset of $\pi$ defined by $\cali{K}_{\pi}=\{\pi\} \cup \cali{K}_{\widetilde{\pi}}$.
If
$\g$ is of exceptional type, $\pi\setminus\widetilde{\pi}$ only
consists of one simple root which we denote by $\alpha_{\pi}$.
Note that $\alpha_{\pi}$ is the simple root which is connected to the lowest root
in the extended Dynkin diagram.

As a consequence of Theorem~\ref{t-tran}, to describe all the
quasi-reductive parabolic subalgebras of $\g$, for $\g$ of exceptional type,
it suffices to consider the case of parabolic subalgebras $\p_{\pi'}^+$ with
$\alpha_{\pi}\in\pi'$.
This will be an important reduction in the sequel.

\begin{remark}\label{r-tran}
If $\g$ has type
$\mathrm{F}_4$ (resp.~$\mathrm{E}_6$, $\mathrm{E}_7$, $\mathrm{E}_8$),
then $\g_{\widetilde{\pi}}$ has type $\mathrm{C}_3$ (resp.~$\mathrm{A}_5$,
$\mathrm{D}_6$, $\mathrm{E}_7$).
In particular, if $\g$ has type $\mathrm{F}_{4}$ or $\mathrm{E}_6$,
then $\p_{\pi'}^+$ is quasi-reductive for any $\pi'$ which does not contain $\alpha_{\pi}$ because
in types A and C all (bi)parabolic subalgebras are quasi-reductive.
\end{remark}

%
\subsection{}
%

As a next step we now focus on a property that we call ``additivity'' to relate the
quasi-reductivity of different parabolic subalgebras (cf.~Theorem~\ref{t-add}).
Throughout this paragraph, $\g$ is assumed to be simple.

\begin{definition}
Let $\pi',\pi''$ be subsets of $\pi$.
We say that $\pi'$ is {\em not connected to} $\pi''$ if $\alpha'$ is orthogonal to $\alpha''$,
for all $(\alpha',\alpha'')$ in $\pi'\times\pi''$.
\end{definition}

\begin{notation}\label{n-K+}
For a positive root $\alpha$,
we denote by $K_{\pi}^+(\alpha)$
the only element $L$ of $\cali{K}_{\pi}$
such that $\alpha\in\Gamma_L$.
Note that unless $\alpha \in \mathcal{E}_{\pi}$,
$K_{\pi}^+(\alpha)$ is the only element $L$ of $\cali{K}_{\pi}$ for which
$\eps{L}-\alpha$ is a positive root.
For $K \in \cali{K}_{\pi}$, we have $K_{\pi}^+(\eps{K})=K$.
\end{notation}

\begin{remark}
It can be checked that $K_{\pi}^+(\alpha)=K_{\pi}^+(\beta)$ for
$\alpha$, $\beta$ simple if and only if $\alpha$ and $\beta$ are in the same orbit of $-w_0$ where $w_0$ is the longest element of
the Weyl group of $\g$.
This suggests that $w_0$ should play a role in these questions, as may be guessed
from a result of Kostant which says that $\mathcal{E}_{\pi}$ is a basis of the space of fixed points of $-w_0$ and from work of Joseph and collaborators (\cite{Jo1,Jo2}).
\end{remark}

\begin{definition}
We shall say that two subsets $\pi',\pi''$
which are not connected to each other \emph{satisfy the condition $(*)$} if:
\begin{eqnarray*}
(*)
&&
K_{\pi}^{+}(\alpha')\not= K_{\pi}^{+}(\alpha'')\quad \forall\ (\alpha',\alpha'')\in \pi' \times \pi''\, .
\end{eqnarray*}
\end{definition}
Note that if $\kg{\pi}=\rk\g$ (that is if $-w_0$ acts trivially on $\pi$), the condition
$(\ast)$ is always satisfied.
Moreover, by using Table \ref{T-exc}, a case-by-case discussion shows:
\begin{lemma}\label{l-star}
Assume that $\g$ is simple of exceptional type and let $\pi'$ be a connected
subset of $\pi$ containing $\alpha_{\pi}$.
Then, for any subset $\pi''$ of $\pi$ which is not connected to $\pi'$,
the two subsets $\pi'$, $\pi''$ satisfy the condition $(*)$, unless $\g=\mathrm{E}_6$,
$\pi'=\{\alpha_1,\alpha_2,\alpha_3,\alpha_4\}$ and $\pi''=\{\alpha_6\}$
or by symmetry $\pi'=\{\alpha_2,\alpha_4,\alpha_5,\alpha_6\}$ and $\pi''=\{\alpha_1\}$.
\end{lemma}

\begin{remark}
If $\g=\mathrm{E}_6$, with
$\pi'=\{\alpha_1,\alpha_2,\alpha_3,\alpha_4\}$ and $\pi''=\{\alpha_6\}$, then
$K_{\pi}^+(\alpha_1)=K_{\pi}^+(\alpha_6)=\{\{\alpha_1,\alpha_3,\alpha_4,\alpha_5,\alpha_6\}\}$, so
$\pi'$ and $\pi''$ do not satisfy the condition $(\ast)$.
As a matter of fact, the parabolic subalgebra $\p_{\pi'\cup\pi''}^+$
will appear as a very special case (see Remark~\ref{r-add}).
\end{remark}

Let $\pi',\pi''$ be two subsets of $\pi$ which are not connected to each other
and assume that $\pi',\pi''$ satisfy condition $(\ast)$.
By Proposition~\ref{p-bor}(ii), we can let $w'$ be in $\l_{\pi'}$ such that $\varp{w}{\p_{\pi'}^+}$ is regular
where $w=w'+\uq{\pi'}$.
Denote by $\s'$ be the image of $\p_{\pi'}^+(\var{w})$ by the projection map from
$\p_{\pi'}^+$ to its derived Lie algebra $\g_{\pi'}\oplus\m_{\pi'}^+$
with respect to the decomposition $\p_{\pi'}^+=\z(\l_{\pi'})\oplus\g_{\pi'}\oplus\m_{\pi'}^+$.
Let $\mathfrak{k}'$ be the intersection of
$\z(\l_{\pi'})$ with $\bigcap\limits_{\varepsilon\in\cali{E}_{\pi},\,\eps{}\not\in\Delta_{\pi'}^+}
\ker \eps{}$.

\begin{lemma}\label{l-add}
{\rm (i)} $\ind \p_{\pi'}^+=\dim \s' +\dim \mathfrak{k}'$.

{\rm (ii)} $[\s', \p_{\pi'\cup \pi''}^+] \subset \p_{\pi'}^+$ and
$\var{w}([\s', \p_{\pi' \cup \pi''}^+])=\{0\}$.
\end{lemma}

\begin{proof}
(i) We have $\dim \p_{\pi'}^+(\var{w})=\ind \p_{\pi'}^+$.
Since the image of $\p_{\pi'}^+(\var{w})$ by the projection map
from $\p_{\pi'}^+$ to $\g_{\pi'}\oplus\m_{\pi'}^+$ is $\mathfrak{s}'$,
it suffices to observe that the intersection of
$\z(\l_{\pi'})$ with $\p_{\pi'}^+(\var{w})$ is $\mathfrak{k}'$.
And this follows from the choice of $w$.

(ii) Let $x$ be an element of $\p_{\pi'}^+(\var{w})$;
write $x=x_0+x'+x^+$ with $x_0\in \z(\l_{\pi'})$, $x'\in\g_{\pi'}$ and $x^+ \in \m_{\pi'}^+$.
Since $[x^+,w']$ lies in $\m_{\pi'}^+$, the fact that $x \in \p_{\pi'}^+(\var{w})$ means
$[x_0,\uq{\pi'}]+[x',w'] +[x',\uq{\pi'}] +[x^+,\uq{\pi'}] \in \m_{\pi'}^+$.
First, we have to show $[x'+x^+, \p_{\pi'\cup \pi''}^+] \subset \p_{\pi'}^+$.
As $[x',\p_{\pi'\cup \pi''}^+]  \subset \p_{\pi'}^+$ since $\pi'',\pi'$ are not connected,
it suffices to prove that $x^{+} \in \m_{\pi'\cup \pi''}^+$.
If not, there are $\gamma \in \Delta_{\pi''}^+$, $K \in \cali{K}_{\pi}$, and $\alpha' \in \Delta_{\pi'}^+$  such that
$$\gamma - \eps{K_{\pi}^{+}(\gamma)} = - (\alpha' + \eps{K})  \mbox{ ,}
\textrm{ i.e. }
\eps{K_{\pi}^{+}(\gamma)} = \gamma + (\alpha' + \eps{K})$$
Hence $\gamma, \alpha' \in \Gamma_{K_{\pi}^{+}(\gamma)}^{0}$ that is
$K_{\pi}^{+}(\alpha')=K_{\pi}^{+}(\gamma)$.
But this contradicts condition $(\ast)$.
Thus $[x'+x^+, \p_{\pi'\cup \pi''}^+]  \subset \p_{\pi'}^+$.

It remains to show: $\var{w}([x'+x^+, \p_{\pi' \cup \pi''}^+])=\{0\}$
that is $[x'+x^+,w] \in \m_{\pi'\cup \pi''}^+$.
If $[x'+x^+,w] \not\in \m_{\pi'\cup \pi''}^+$,
there must be $\gamma \in \Delta_{\pi}^+ \setminus \Delta_{\pi'}$,
$K \in \cali{K}_{\pi}$, and $\alpha'' \in \Delta_{\pi''}^+$
such that $\gamma - \eps{K} = \alpha''$.
In particular $\alpha'' \in \Gamma_{K_{\pi}^{+}(\gamma)}^{0}$
that is $K_{\pi}^{+}(\alpha'')=K_{\pi}^{+}(\gamma)$.
On the other hand, $[x,w] \in \m_{\pi'}^+$ implies
that there exist $\alpha' \in  \Delta_{\pi'}^+$ and $L \in \cali{K}_{\pi}$,
such that
$$\gamma- \eps{K_{\pi}^{+}(\gamma)} = - (\alpha' + \eps{L}) \mbox{ ,}
\textrm{ i.e. }
\eps{K_{\pi}^{+}(\gamma)} = \gamma + (\alpha' + \eps{L})$$
As before, we deduce that $\alpha' \in \Gamma_{K_{\pi}^{+}(\gamma)}^{0}$,
i.e.~$K_{\pi}^{+}(\alpha')=K_{\pi}^{+}(\gamma)=K_{\pi}^{+}(\alpha'')$
and this contradicts condition $(\ast)$.
\end{proof}

\begin{corollary} \label{c-add}
Let $\pi',\pi''$ be two subsets of $\pi$ which are not connected to each other
and satisfy condition $(\ast)$.
If $\p_{\pi'\cup \pi''}^+$ is quasi-reductive then $\p_{\pi'}^+$ and $\p_{\pi''}^+$ are both quasi-reductive.
\end{corollary}

\begin{proof}
Suppose that $\p_{\pi'\cup \pi''}^+$ is quasi-reductive and
that any one of the other two parabolic subalgebras is not quasi-reductive
and show that this leads to a contradiction.
By assumption we can choose $\varphi\in (\p_{\pi'\cup \pi''}^+)^*$ of reductive type for $\p_{\pi'\cup \pi''}^+$
such that $\varphi'=\varphi_{|_{\p_{\pi'}^+}}$ and $\varphi''=\varphi_{|_{\p_{\pi''}^+}}$ are $\p_{\pi'}^+$-regular
and  $\p_{\pi''}^+$-regular respectively.
Suppose for instance that $\p_{\pi'}^+$ is not quasi-reductive.
By Proposition~\ref{p-bor}(ii) we can suppose furthermore
that $\varphi'=\varp{w}{\p_{\pi'}^+}$ for some $w=w'+\uq{\pi'}$ with $w'\in \l_{\pi'}$.

Since we assumed that $\p_{\pi'}^+$ is not quasi-reductive, $\varp{w}{\p_{\pi'}^+}$ contains a nonzero nilpotent element, $x$,
which is so contained in the derived Lie algebra of $\p_{\pi'}^+$.
Then, Lemma~\ref{l-add}(ii) gives $[x, \p_{\pi'\cup \pi''}^+] \subset \p_{\pi'}^+$ and
$\{0\}=\var{w}([x, \p_{\pi'\cup \pi''}^+])= \varphi'([x, \p_{\pi'\cup \pi''}^+])=\varphi([x, \p_{\pi'\cup \pi''}^+])$.
As a consequence, $\p_{\pi'\cup\pi''}^+(\varphi)$ contains the nonzero nilpotent element $x$.
This contradicts the choice of $\varphi$.
The same line of arguments works if we assume that $\p_{\pi''}^+$ is not quasi-reductive.
\end{proof}

Under certain conditions, the converse of Corollary~\ref{c-add} is also true as we show now.
To begin with, let us express  the index of $\p_{\pi'\cup \pi''}^+$ in terms of those of $\p_{\pi'}^+$ and $\p_{\pi''}^+$.
As $E_{\pi'\cup \pi'',\pi}=E_{\pi',\pi}+E_{\pi',\pi}$, we get:
$\dim E_{\pi'\cup \pi'',\pi} = \dim E_{\pi',\pi} + \dim E_{\pi'',\pi} - \dim (E_{\pi',\pi} \cap E_{\pi'',\pi})$.
Hence, formula~(\ref{eq-ind}) implies
\begin{eqnarray}\label{e-add}
\ind \p_{\pi'\cup \pi''}^+&=& \ind \p_{\pi'}^+ +\ind \p_{\pi''}^+ - 
\left(\rk\g +\kg{\pi} -2 \dim (E_{\pi',\pi} \cap E_{\pi'',\pi} ) \right) \mbox{ . }
\end{eqnarray}
In case $\rk\g=\kg{\pi}$, the intersection $E_{\pi',\pi} \cap E_{\pi'',\pi}$ is equal to $E_{\pi}$ and has dimension $\rk\g$.
Hence, the index is additive in that case, as~(\ref{e-add}) shows.\\

\begin{theorem}[Additivity] \label{t-add}
Assume that $\g$ is simple and of exceptional type and that $\kg{\pi} =\rk\g$.
Let $\pi',\pi''$ be two subsets of $\pi$
which are not connected to each other.
Then, $\p_{\pi'\cup \pi''}^+$  is quasi-reductive
if and only if both $\p_{\pi'}^+$ and $\p_{\pi''}^+$ are quasi-reductive.
\end{theorem}

\begin{remark} \label{r-add}
The conclusions of Theorem~\ref{t-add} is valid for classical simple Lie algebras, even without the hypothesis
$\kg{\pi}=\rk\g$.
In types $\mathrm{A}$ or $\mathrm{C}$ this follows from the fact that all biparabolic subalgebras are quasi-reductive.
If $\g$ is an orthogonal Lie algebra, this is a consequence of Theorem~\ref{t-dkt}.
However, for the exceptional Lie simple algebra E$_6$, the only one for which $\kg{\pi}\not=\rk\g$,
the conclusions of Theorem~\ref{t-add} may fail.
Indeed, let us consider the following subsets of $\pi$ for $\g$ of type E$_6$:
$\pi'=\{\alpha_1,\alpha_2,\alpha_3,\alpha_4\}$ and $\pi''=\{\alpha_6\}$..
By Remark~\ref{r-0}, $\p_{\pi'}^+$ is quasi-reductive as a Lie algebra of zero index.
On the other hand, $\p_{\pi''}^{+}$ is quasi-reductive by the transitivity property, cf.~Remark~\ref{r-tran}.
But, it will be shown in Theorem~\ref{t-nE6} that $\p_{\pi' \cup\pi''}^{+}$ is not quasi-reductive.

As a consequence of Lemma~\ref{l-star} and Corollary~\ref{c-add},
even in type E$_6$ where $\rk\g\not=\kg{\pi}$, if $\p_{\pi'\cup\pi''}^+$ is quasi-reductive, then $\p_{\pi'}^+$ and $\p_{\pi''}^+$ are both quasi-reductive.

As a by-product of our classification, we will see that the above situation is the only case which prevents the additivity property to be true for all simple Lie algebras (see Remark~\ref{r-addE6}).
\end{remark}

\begin{proof}
We argue by induction on the rank of $\g$.
By the transitivity property (Theorem~\ref{t-tran}), Remark~\ref{r-add} and the induction, we can assume that $\alpha_{\pi} \in\pi'$.
Then, by Lemma~\ref{l-star} and Corollary~\ref{c-add}, only remains to prove that
if both $\p_{\pi'}^+$ and $\p_{\pi''}^+$ are quasi-reductive, then so is $\p_{\pi'\cup \pi''}^+$.

Assume that both $\p_{\pi'}^+$ and $\p_{\pi''}^+$ are quasi-reductive.
By Proposition~\ref{p-Zar}, we can find a linear regular form $\varphi$ in $(\p_{\pi'\cup \pi''}^+)^*$
such that $\varphi'=\varphi_{|_{\p_{\pi'}^+}}$ and $\varphi''=\varphi_{|_{\p_{\pi''}^+}}$ are regular and of reductive type for
$\p_{\pi'}^+$ and $\p_{\pi''}^+$ respectively.
By Proposition \ref{p-bor}(ii), we can assume that $\varphi=\varp{w+\uq{}}{\p_{\pi'\cup \pi''}^+}$,
where $w=h+w'+w''$, with $w'\in \n_{\pi'}^+$, $w''\in \n_{\pi''}^+$ and $h \in\h$.
Hence, $\varphi'=\varp{h + w'+ \uq{}}{\p_{\pi'}^+}$ and $\varphi''=\varp{h + w''+ \uq{}}{\p_{\pi''}^+}$.

Use the notations of Lemma~\ref{l-add}.
By Lemma~\ref{l-add}(ii), $\s'$ is contained in $\p_{\pi'\cup\pi''}^+(\varphi)$.
Show now that $\mathfrak{k}'$ is zero.
Let $h$ be an element $\mathfrak{k}'$.
Since $h\in \mathfrak{k}'$, we have $\eps{}(h)=0$ for any $\eps{}\in \cali{E}_{\pi}$
which is not in $\Delta_{\pi'}^+$.
On the other hand, for any $\eps{} \in \cali{E}_{\pi} \cap \Delta_{\pi'}^+$,
we have $\eps{}(h)=0$ since $h$ lies in the center of $\l_{\pi'}$.
Hence, our assumption $\rk\g=\kg{\pi}$ implies $h=0$.
As a consequence of Lemma~\ref{l-add}(i), we deduce that $\ind\p_{\pi'}^+=\dim\s'$.
Similarly, if $\s''$ denotes the image of $\p_{\pi''}^+(\varphi'')$ under the projection from $\p_{\pi''}^+$ to $\g_{\pi''}\oplus\m_{\pi''}^+$,
Lemma~\ref{l-add}(ii) tells us that $\s''$ is contained in $\p_{\pi'\cup\pi''}^+(\varphi)$ and that $\ind\p_{\pi''}^+=\dim\s''$.

To summarize, our discussion shows that $\s' + \s''$ is contained in $\p_{\pi'\cup\pi''}^+(\varphi)$
and that these two subspaces have the same dimension by equation~(\ref{e-add}).
So $\s' + \s'' = \p_{\pi'\cup\pi''}^+(\varphi)$.
But by assumption, $\s' + \s'' $ only consists of semisimple elements.
From that we deduce that $\varphi$ is of reductive type for $\p_{\pi'\cup\pi''}^+$, whence the theorem.
\end{proof}

\section{Some classes of quasi-reductive biparabolic subalgebras}\label{S-3}

In this section we show that, under certain conditions on the interlacement of the two cascades
of $\pi_1$ and $\pi_2$, we can deduce that $\q_{\pi_1,\pi_2}$ is quasi-reductive
(Theorem~\ref{t-lac}).
We assume in this section that $\g$ is simple.

%
\subsection{}
%

We start by introducing the necessary notations.
Recall that for a positive root $\alpha$, $K_{\pi}^+(\alpha)$ stands for the only element $L$ of $\cali{K}_{\pi}$
such that $\alpha\in\Gamma_L$, cf.~Notation~\ref{n-K+}.
To any positive root $\alpha\in \Delta_{\pi}^+$
we now associate the subset $\cali{K}_{\pi}^-(\alpha)$ of the cascade $\mathcal{K}_{\pi}$ of all $L$
such that the highest root $\eps{L}$ can be added to $\alpha$:
\begin{eqnarray*}
\cali{K}_{\pi}^-(\alpha) &=&\{L\in \cali{K}_{\pi} \, | \, \eps{L}+\alpha \in  \Delta_{\pi}^{+}\}\, .
\end{eqnarray*}
Observe that the set $\cali{K}_{\pi}^-(\alpha)$ may be empty or contain more than one element.

\begin{examples}
(1) If $K$ is in the cascade $\cali{K}_{\pi}$ then $\cali{K}_{\pi}^-(\eps{K})$ is empty.

(2) In type $\mathrm{E}_7$, for $\alpha=\alpha_4 +\alpha_5 +\alpha_6$, the set $\cali{K}_{\pi}^-(\alpha)$ has more than one element:
$\eps{4}+\alpha,\eps{5}+\alpha,\eps{6}+\alpha$ are all positive roots.
\end{examples}

\noindent
We need also the following notation:
$$\widetilde{\Delta}^+_{\pi}=\{\alpha \in \Delta^{+}_{\pi}, \, \alpha = \frac{1}{2}(\eps{K}-\eps{K'}) \, ; \,  K,K' \in \cali{K}_{\pi}\}\, .$$

\begin{remark} \label{r-s-l}
One can check that for $\g$ a simply-laced simple Lie algebra,
no positive root can be written in the way as asked for in the definition of $\widetilde{\Delta}^+_{\pi}$.
Thus $\widetilde{\Delta}^+_{\pi}$ is empty if $\g$ is simple of type $\mathrm{A}$, $\mathrm{D}$ or $\mathrm{E}$.
\end{remark}

We list the sets $\widetilde{\Delta}^+_{\pi}$ in Table~\ref{t-eps} for
the simple Lie algebras of types $\mathrm{B}_{\ell}$, $\mathrm{C}_{\ell}$,
$\mathrm{G}_2$ and $\mathrm{F}_4$.

{\begin{table}[r,h] \tiny
\begin{center}
\begin{tabular}{ll}
\hline
\\
$\mathrm{B}_{\ell}$, $\ell \geq 2$:  & $\{\frac{1}{2}(\varepsilon_{2i}-\varepsilon_{2i-1})$, $i=1,\ \ldots,\, \left[\frac{\ell}{2}\right]\}$ \\

\medskip

$\mathrm{C}_{\ell}$, $\ell \geq 3$:  & $\{\frac{1}{2}(\varepsilon_{i} - \varepsilon_{i+k+1})$,  $1 \leq i \leq  \ell-1$, $0\leq k \leq l-i-1\}$\\

\medskip

$\mathrm{G}_2$:  & $\{\begin{array}{c}11\end{array}=\frac{1}{2}(\varepsilon_{1} - \varepsilon_{2})\}$ \\

\medskip
$\mathrm{F}_4$: &
$\{\begin{array}{c}1110\end{array}=\frac{1}{2}(\varepsilon_{1} - \varepsilon_{2})$,
$\begin{array}{c}1111\end{array}=\frac{1}{2}(\varepsilon_{1} - \varepsilon_{3})$,
$\begin{array}{c}1121\end{array}=\frac{1}{2}(\varepsilon_{1} - \varepsilon_{4})$,\\
& \hspace{0.05cm} $\begin{array}{c}0001\end{array}=\frac{1}{2}(\varepsilon_{2} - \varepsilon_{3})$,
$\begin{array}{c}0011\end{array}=\frac{1}{2}(\varepsilon_{2} - \varepsilon_{4})$,
$\begin{array}{c}0010\end{array}=\frac{1}{2}(\varepsilon_{3} - \varepsilon_{4})\}$\\
\\
\hline
\end{tabular}

\vspace{.02cm}
\caption{\label{t-eps} $\widetilde{\Delta}^+_{\pi}$ for the simple Lie algebras.}
\end{center}
\end{table}}

Part of the following lemma explains that for a root $\alpha$ in $\widetilde{\Delta}^+_{\pi}$
we can actually describe the two cascades involved in the expression of $\alpha$:

\begin{lemma}\label{l-cas}
{\rm (i)} Whenever $\alpha \in \widetilde{\Delta}^{+}_{\pi}$,
then $\cali{K}_{\pi}^-(\alpha)$ consists of a unique element $K_{\pi}^-(\alpha)$.

{\rm (ii)} For any element $\alpha= \frac{1}{2}(\eps{K}-\eps{K'})$ of $\widetilde{\Delta}^+_{\pi}$
we have $K=K_{\pi}^+(\alpha)$ and $K'=K_{\pi}^-(\alpha)$.
\end{lemma}

\begin{proof}
One can deduce (i) from Table~\ref{t-eps}.

(ii) By (i), we have $\cali{K}_{\pi}^-(\alpha)=\{K_{\pi}^-(\alpha)\}$.
Furthermore, \hbox{$<\alpha,\eps{K}^{\vee}>=1$} so $\eps{K}-\alpha$ is a root (cf.~\cite[Proposition 18.5.3(iii)]{TY1}).
Since $\eps{K}-\alpha=\eps{K'}+\alpha$, these two are both positive roots, forcing $K_{\pi}^+(\alpha)=K$ and $K_{\pi}^-(\alpha)=K'$.
\end{proof}

Let $\pi_1$ and $\pi_2$ be two subsets of $\pi$.
We define
$$\widetilde{\cali{K}}_i^{(j)}=\{M\in\cali{K}_{\pi_i} \ | \ \eps{M}\in\widetilde{\Delta}^+_{\pi_j}\}\, .$$
Thus, for $M$ is in $\widetilde{\cali{K}}_{i}^{(j)}$ we have
$\eps{M}=\frac{1}{2}(\eps{K_{\pi_j}^{+}(\eps{M})} - \eps{K_{\pi_j}^{-}(\eps{M})})$
by Lemma~\ref{l-cas}(ii).
Note that $M$ is an element of the cascade of $\pi_i$ while $K_{\pi_j}^{\pm}(\eps{M})$ belong to the cascade of $\pi_j$.

\begin{definition}\label{d-lac}
Let $\pi_1$, $\pi_2$ be subsets of $\pi$.
We say that the cascades $\cali{K}_{\pi_1}$ and $\cali{K}_{\pi_2}$ are {\em well-interlaced}
if $\dim (E_{\pi_1} \cap E_{\pi_2})
=\#(\cali{K}_{\pi_1} \cap \cali{K}_{\pi_2}) + \# \widetilde{\cali{K}}_{1}^{(2)}+  \# \widetilde{\cali{K}}_{2}^{(1)}$.
\end{definition}

\begin{remark} \label{r-cas}
The following subsets $\pi_1,\pi_2$ of $\pi$ give rise to examples of well-interlaced cascades:

(1) $\pi_1$ and $\pi_2$ are such that $\cali{K}_{\pi_i} \subset \cali{K}_{\pi_j}$ or $\cali{K}_{\pi_j}\subset \cali{K}_{\pi_i}$.
In particular, this is the case if $\pi_1$ or $\pi_2$ is empty.

(2) $\pi_1$ and $\pi_2$ are such that the collection
of all highest roots $\cali{E}_{\pi_1} \cup \cali{E}_{\pi_2}$
consists of linearly independent elements\footnote{We mean that
this collection of roots forms
a set of linearly independent roots, neglecting
any multiplicities that might occur, cf.~Example~\ref{ex2-E6} below.}.

These two cases have already been studied by Tauvel and Yu in~\cite{TY2}.
\end{remark}

We are now ready to formulate the main result of this section. It will
be proved in Subsection~\ref{Ss-proof} below.

\begin{theorem}\label{t-lac}
Let $\q_{\pi_1,\pi_2}$ be a biparabolic subalgebra of $\g$.
Assume that the cascades $\cali{K}_{\pi_1}$ and $\cali{K}_{\pi_2}$ are well-interlaced.
Then $\mathfrak{q}_{\pi_1,\pi_2}$ is quasi-reductive.

More precisely, the linear form $\varphi_{u(\underline{a},\underline{b})}$ is of reductive type
for almost all choices of the coefficients $(\underline{a}, \underline{b})\in\C^{\kg{\pi_1}+\kg{\pi_2}}$.
\end{theorem}

\begin{example} \label{ex2-E6}
Suppose that $\g$ is simple of type $\mathrm{E}_6$.
In the case where $\pi_1=\{\alpha_2,\alpha_3,\alpha_4\}$
(resp.~$\pi_1=\{\alpha_2,\alpha_3,\alpha_4,\alpha_6\}$, $\pi_1=\{\alpha_1,\alpha_2,\alpha_3,\alpha_4\}$) and $\pi_2=\pi$,
the union $\cali{E}_{\pi_1} \cup \cali{E}_{\pi_2}$ consists of linearly independent elements.
Hence $\q_{\pi_1,\pi_2}=\p_{\pi_1}^{+}$ is quasi-reductive by Remark~\ref{r-cas}(2) and Theorem~\ref{t-lac}.
\end{example}

We now give an example which is not covered by 
Remark~\ref{r-cas}:

\begin{example}
Suppose that $\g$ is simple of type $\mathrm{F}_4$.
The subsets $\pi_1=\{\alpha_3,\alpha_4\}$ and $\pi_2=\pi$ are well-interlaced
and $\q_{\pi_1,\pi_2}=\p_{\pi_1}^{+}$ is quasi-reductive by Theorem~\ref{t-lac}.
Note that Theorem~\ref{t-tran} provides an alternative way to prove that this parabolic subalgebra is quasi-reductive.
\end{example}

\begin{remark}
The converse of Theorem~\ref{t-lac} is not true.
For example, we can easily check that the assumption of Theorem \ref{t-lac}
does not hold for the parabolic subalgebra $\p_{\{\alpha_2,\alpha_4\}}^+$ of $\mathrm{E}_6$.
However, it is quasi-reductive as we will show in Subsection~\ref{Ss-2} (Theorem~\ref{t-2}).
\end{remark}

%
\subsection{}\label{Ss-proof}
%

This subsection is devoted to the proof of Theorem~\ref{t-lac}.
We start with two technical lemmata.

Let $\alpha \in \widetilde{\Delta}_{\pi}^{+}$.
Recall that by Lemma~\ref{l-cas}(ii), $\alpha$ is written as
$\alpha=\frac{1}{2}(\eps{K_{\pi}^+(\alpha)}-\eps{K_{\pi}^-(\alpha)})$.
As an abbreviation we set
$$\overline{\alpha}=\frac{1}{2}(\eps{K_{\pi}^+(\alpha)}+\eps{K_{\pi}^-(\alpha)})\, .$$
Thus $\alpha +\overline{\alpha}=\eps{K_{\pi}^{+}(\alpha)}$ and $-\alpha+\overline{\alpha}=\eps{K_{\pi}^-(\alpha)}$.
From the relations between the four roots $\alpha$, $\overline{\alpha}$, $\eps{K_{\pi}^-(\alpha)}$ and $\eps{K_{\pi}^+(\alpha)}$
we define the structure constants $\tau_1, \tau_2,\tau_3,\tau_4$ as follows:

\begin{center}
\begin{tabular}{ll}
$[ x_{\alpha}, x_{-\eps{K_{\pi}^+(\alpha)}}  ]=\tau_1  x_{-\overline{\alpha}}$ \ ; &
$[ x_{-\alpha}, x_{-\eps{K_{\pi}^-(\alpha)}}  ]=\tau_2  x_{-\overline{\alpha}}$\ ; \\
$[  x_{\alpha}, x_{\overline{\alpha}}  ] = \tau_3  x_{\eps{K_{\pi}^{+}(\alpha)}}$ \ ; &
$[ x_{-\alpha},x_{\overline{\alpha}} ] = \tau_4  x_{ \eps{K_{\pi}^{-}(\alpha)}}$ \ .
\end{tabular}
\end{center}

\begin{lemma}\label{l-tec1}
Assume that $\g$ is of type $\mathrm{B}_{\ell}$ ($\ell\geq 2$), $\mathrm{C}_{\ell}$ ($\ell\geq 3$) or $\mathrm{F}_4$.
Let $\alpha$ be in $\widetilde{\Delta}^{+}_{\pi}$.

{\rm (i)} The only roots of the form $k\alpha+l\overline{\alpha}$ are $\{\pm\alpha,\ \pm\overline{\alpha},\ \pm(\alpha\pm \overline{\alpha})\}$.

{\rm (ii)} We have $\tau_1,\tau_2 \in \{-1,1\}$, $\tau_3,\tau_4 \in \{-2,2\}$ and $\tau_1\tau_4=\tau_2\tau_3$.
\end{lemma}

\begin{proof}
(i) By assumption, the four linear
combinations $\pm(\alpha\pm \overline{\alpha})$ are all roots.
The claim then follows
since root strings have at most length $2$ in types $\mathrm{B}$,
$\mathrm{C}$ and $\mathrm{F}$.

(ii) We explain how to obtain $\tau_1=\pm 1$, the computations of
$\tau_i$ for $i=2,3,4$ is completely analogous.
Consider the $\alpha$-string through $-\eps{K_{\pi}^+(\alpha)}$.
It has the form $\{-\eps{K_{\pi}^+(\alpha)},-\overline{\alpha}, -\eps{K_{\pi}^-(\alpha)}\}$,
so in particular, $p=0$ in the notation of Subsection~\ref{Ss-1.2}, whence $\tau_1=\pm 1$.

Only remains to proof the equality $\tau_1\tau_4=\tau_2\tau_3$.
We compute the bracket $[x_{-\alpha},[x_{\alpha},x_{\overline{\alpha}}]]$ in two different ways.
We have $[x_{-\alpha}, x_{\eps{K_{\pi}^{+}(\alpha)}} ]=-\tau_1 x_{\overline{\alpha}}$ (cf.~\cite[\S18.2.2 and Corollary 18.5.5]{TY1}).
Hence $[x_{-\alpha},[x_{\alpha},x_{\overline{\alpha}}]] = [x_{-\alpha}, \tau_3  x_{\eps{K_{\pi}^{+}(\alpha)}} ]
= -\tau_1 \tau_3  x_{\overline{\alpha}}$.
On the other hand, as $\eps{K_{\pi}^{+}}$ and $\eps{K_{\pi}^{-}}$ have the same length
($\g$ having type different from G$_2$), we have:
$[h_{\alpha},x_{\overline{\alpha}}]=\langle \overline{\alpha} ,\alpha^{\vee} \rangle  x_{\overline{\alpha}} = 0$.
So: $[x_{-\alpha},[x_{\alpha},x_{\overline{\alpha}}]] = [x_{\overline{\alpha}},[x_{\alpha},x_{-\alpha}]] + [x_{\alpha},[x_{-\alpha},x_{\overline{\alpha}}]]
= - [h_{\alpha},x_{\overline{\alpha}}]   +\tau_4  [x_{\alpha}, x_{\eps{K_{\pi}^{-}(\alpha)}}] = - \tau_4\tau_2  x_{\overline{\alpha}}$
again by using~\cite[\S18.2.2 and Corollary 18.5.5]{TY1}.
We have so obtained $\tau_1\tau_3 =\tau_2\tau_4$.
From that the claim follows.
\end{proof}

From now, we let $\pi_1,\pi_2$ be two subsets of $\pi$.

\begin{lemma}\label{l-tec2}
Let $M$ be an element of $\widetilde{\cali{K}}_{i}^{(j)}$.

{\rm  (i)} $\eps{K_{\pi_j}^{+}(\eps{M})}$ and $\eps{K_{\pi_j}^{-}(\eps{M})}$ are not roots of $\pi_i$.

{\rm (ii)} For $K \in \cali{K}_{\pi_j}$, $\eps{M} \pm \eps{K}$ is a root if and only if $K=K_{\pi_j}^{\mp}(\eps{M})$.
\end{lemma}

\begin{proof}
(i) Can be deduced from Tables~\ref{T-cla},~\ref{T-exc} and~\ref{t-eps}.

(ii) The fact that $\eps{M} - \eps{K_{j}^{+}(\eps{M})}$ and $\eps{M} + \eps{K_{j}^{-}(\eps{M})}$ are roots of $\pi_j$
has been observed in Lemma \ref{l-tec1}(i).
Next, by Lemma~\ref{l-cas}(i), we know that $K_{\pi_j}^{-}(\eps{M})$ is the only element $L$ of $\cali{K}_{\pi_j}$
such that $\eps{M}+\eps{L}$ is a root.
Suppose now that there is $L \in \cali{K}_{\pi_j}$, $L \not= K_{\pi_j}^+(\eps{M})$, such that $\eps{L}-\eps{M}$ is a root.
By Lemma \ref{l-tec1}(i), we have $L \not=K_{\pi_j}^{-}(\eps{M})$.
So, the fact that $\eps{L}-\eps{M}$ is a root forces $\beta=\eps{M} - \eps{L}$ to be a positive root,
by definition of $K_{\pi_j}^+(\eps{M})$.
Then the equality $\beta + \eps{L}=\eps{M}$ implies $\langle \beta,\eps{M}^{\vee} \rangle=1$.
On the other hand, we have
$\langle \eps{M},\eps{L}^{\vee}\rangle
= \langle \frac{1}{2} (\eps{K_{\pi_j}^{+}(\eps{M})}-\eps{K_{\pi_j}^{-}(\eps{M})}), \eps{L}^{\vee}\rangle= 0$
since $L \not= K_{\pi_j}^{\pm}(\eps{M})$.
So $\langle \eps{L},\eps{M}^{\vee}\rangle=0$.
As a consequence,
$1=\langle \beta,\eps{M}^{\vee}\rangle= \langle \beta+\eps{L},\eps{M}^{\vee}\rangle =\langle \eps{M},\eps{M}^{\vee}\rangle=2$.
Hence we get a contradiction.
\end{proof}

Recall that for $(\underline{a},\underline{b})\in (\C^*)^{\kg{\pi_1}+\kg{\pi_2}}$, we have set
$$u(\underline{a},\underline{b})=\sum\limits_{K \in \mathcal{K}_{\pi_2}} a_K x_{-\eps{K}}
+ \sum\limits_{L \in \mathcal{K}_{\pi_1}} b_L x_{\eps{L}}$$

\begin{lemma} \label{l-yzt}
Let $(\underline{a},\underline{b})$ be in $(\C^*)^{\kg{\pi_1}+\kg{\pi_2}}$.
For $K \in \cali{K}_{\pi_i} \cap \cali{K}_{\pi_j}$, $M \in  \widetilde{\cali{K}}_1^{(2)}$ and $N \in  \widetilde{\cali{K}}_2^{(1)}$,
there exist $\rho_K \in \C^*$, $(\lambda_M,\mu_M,\nu_M) \in (\C^*)^3$ and $(\lambda'_M,\mu'_M,\nu'_M) \in (\C^*)^3$
such that the elements $y_K$, $z_{M}$ and $t_N$ of $\q_{\pi_1,\pi_2}$ defined by

$\begin{array}{l}
y_K=x_{\eps{K}} + \rho_K  x_{-\eps{K}}\\
z_{M} =  x_{\eps{M}} + \lambda_M  x_{-\eps{M}{}} + \mu_M  x_{\eps{K_{\pi_2}^{+}(\eps{M})}}
+ \nu_M  x_{\eps{K_{\pi_2}^{-}(\eps{M})}}\\
t_{N} = x_{-\eps{N}} + \lambda'_{N}  x_{\eps{N}{}} + \mu'_{N}  x_{-\eps{K_{\pi_1}^{+}(\eps{N})}}
+ \nu'_{N}  x_{-\eps{K_{\pi_1}^{-}(\eps{N})}}
\end{array}$\\
are semisimple elements of $\mathfrak{g}$ which stabilize $\var{u(\underline{a},\underline{b})}$ in $\q_{\pi_1,\pi_2}$.
\end{lemma}

\begin{proof}
Set $u=u(\underline{a},\underline{b})$.
For $K \in \cali{K}_{\pi_i} \cap \cali{K}_{\pi_j}$,
it is clear that $y_K$ is semisimple.
Moreover, for $\rho_K=a_K/b_K$, the element $y_K$ stabilizes $\varp{u}{\q_{\pi_1,\pi_2}}$; we even have $[y_K,u]=0$.

Let now $M$ be in $\widetilde{\cali{K}}_1^{(2)}$.
If $\cali{K}_{1}^{(2)}\ne \emptyset$ then $\g$ cannot be of type G$_2$,
since for G$_2$, $\cali{K}_{1}^{(2)}=\emptyset$ (cf.~Table~\ref{t-eps}).
So $\g$ is of type B$_{\ell}$, C$_{\ell}$ or F$_4$ (Remark~\ref{r-addE6}).
Thus we are in the situation of Lemma~\ref{l-tec1}.
Let $(\lambda_M,\mu_M,\nu_M)$ be in $(\C^*)^3$.
By definition of $z_M$, we have:
$$\begin{array}{l}
[z_M,u] = \sum\limits_{K \in \cali{K}_{\pi_2}} a_K  ([x_{\eps{M}},x_{-\eps{K}}]
+  \lambda_M  [x_{-\eps{M}},x_{-\eps{K}}]) + \mu_M \sum\limits_{L \in \cali{K}_{\pi_1}} b_L [x_{\eps{K_{\pi_2}^{+}(\eps{M})}},x_{\eps{L}}]\\
\quad\quad +  \quad \nu_M \sum\limits_{L \in \cali{K}_{\pi_1}} b_L  [x_{\eps{K_{\pi_2}^{-}(\eps{M})}},x_{\eps{L}}]
-  \lambda_M b_M h_{\eps{M}} + \mu_M a_{K_{\pi_2}^{+}(\eps{M})} h_{\eps{K_{\pi_2}^{+}(\eps{M})}}
+ \nu_M a_{K_{\pi_2}^{-}(\eps{M})}  h_{\eps{K_{\pi_2}^{-}(\eps{M})}}
\end{array}$$
Note that $[x_{\eps{M}},x_{-\eps{K}}]\neq 0$ if and only if $K=K^+_{\pi_2}(\eps{M})$ by Lemma~\ref{l-cas}.
By Lemma~\ref{l-tec2}(i), the element
$v=\mu_M \sum\limits_{L \in \cali{K}_{\pi_1}} b_L [x_{\eps{K_{\pi_2}^{+}(\eps{M})}},x_{\eps{L}}]
+ \nu_M \sum\limits_{L \in \cali{K}_{\pi_1}} b_L [x_{\eps{K_{\pi_2}^{-}(\eps{M})}},x_{\eps{L}}]$
lies in $\mathfrak{u}_{\pi_1,\pi_2}$.
We set $\overline{\varepsilon_{M}}=\frac{1}{2}(\eps{K_{\pi_j}^{+}(\eps{M})}+\eps{K_{\pi_j}^{-}(\eps{M})})$,
and define the structure constants $\tau_1,\tau_2,\tau_3,\tau_4$ for $\alpha=\eps{M}$ and $\overline{\alpha}=\overline{\varepsilon_{M}}$.
Then, by Lemma~\ref{l-cas} we have
$$\begin{array}{rcl}
[z_M,u]&=&( \tau_1  a_{K_{\pi_2}^{+}(\eps{M})}
+\lambda_M  \tau_2  a_{K_{\pi_2}^{-}(\eps{M})})  x_{-\overline{\varepsilon_{M}}} + v\\
&&-  \lambda_M  b_M  h_{\eps{M}} +  \mu_M  a_{K_{\pi_2}^{+}(\eps{M})}  h_{\eps{K_{\pi_2}^{+}(\eps{M})}}
+  \nu_M \ a_{K_{\pi_2}^{-}(\eps{M})}  h_{\eps{K_{\pi_2}^{-}(\eps{M})}}
\end{array}$$
By Remark~\ref{r-s-l}, the elements of $\cali{E}_{\pi_2}$ form a
basis of $\mathfrak{h}_{\pi_2}^*$ since $\cali{K}_{1}^{(2)} \ne \emptyset$.
So, by Lemma~\ref{l-cas}, we can write
$h_{\eps{M}}=c^+  h_{\eps{K_{\pi_2}^{+}(\eps{M})}} - c^-  h_{\eps{K_{\pi_2}^{-}(\eps{M})}}$ with $c^+,c^- \in \C^*$.
Furthermore, $\eps{K_{\pi_2}^{+}(\eps{M})}$ and $\eps{K_{\pi_2}^{-}(\eps{M})}$
have the same length (they are both long roots, cf.~Table \ref{t-eps}).
So, $c^+=c^-$ (cf.~\cite[\S18.3.3]{TY1}).
Hence
$$\begin{array}{rcl}
[z_M,u]&=&  ( \tau_1  a_{K_{\pi_2}^{+}(\eps{M})} + \lambda_M  \tau_2
a_{K_{2}^{-}(\eps{M})})  x_{-\overline{\varepsilon_{M}}}+ v\\
&&+ \ (- c  \lambda_M  b_M + \mu_M  a_{K_{\pi_2}^{+}(\eps{M})})
h_{\eps{K_{\pi_2}^{+}(\eps{M})}}+ (c  \lambda_M \,b_M
+\nu_M  a_{K_{\pi_2}^{-}(\eps{M})} ) h_{\eps{K_{\pi_2}^{-}(\eps{M})}}
\end{array}$$
As a result, if we take for $\lambda_M$, $\lambda_M=- \tau_1 a_{K_{\pi_2}^{+}(\eps{M})}/ (\tau_2  a_{K_{\pi_2}^{-}(\eps{M})})$
and then for $\mu_M$ and $\mu_M$, $\mu_M= c  \lambda_M  b_M/a_{K_{\pi_2}^{+}(\eps{M})}$ and
$\nu_M=- c  \lambda_M  b_M/a_{K_{\pi_2}^{-}(\eps{M})}$ we obtain that $[z_M,u]=v  \in \mathfrak{u}_{\pi_1,\pi_2}$, i.e.~that
$z_M$ stabilizes $\varp{u}{\q_{\pi_1,\pi_2}}$.
In a similar way, one shows that $t_N$ stabilizes $\varp{u}{\q_{\pi_1,\pi_2}}$.

It remains to prove that $z_M$ is semisimple (and that $t_N$ is semisimple but this can be done in a similar way).
By Lemma \ref{l-tec1}(i), we have
\begin{eqnarray*}
\hspace{-1cm}\exp(t  \mathrm{ad}\,x_{\overline{\varepsilon_{M}}})(x_{\eps{M}}+ \lambda_M x_{-\eps{M}})
&=&x_{\eps{M}}+\lambda_M  x_{-\eps{M}}+t  [x_{\overline{\varepsilon_{M}}},x_{\eps{M}}]+ t  \lambda_M  [x_{\overline{\varepsilon_{M}}},x_{-\eps{M}}]\\
&=& x_{\eps{M}}+\lambda_M  x_{-\eps{M}}- t  \tau_3 \, x_{\eps{K_{\pi_2}^{+}(\eps{M})}}- t
\lambda_M \tau_4  x_{\eps{K_{\pi_2}^{-}(\eps{M})}}
\end{eqnarray*}
for any $t\in \C^*$.
By Lemma \ref{l-tec1}(ii), we have
$\tau_1 \tau_4=\tau_2 \tau_3$.
Therefore it is possible to choose $t$ so that both equalities $- t  \tau_3 = c  \lambda_M  b_M/a_{K_{\pi_2}^{+}(\eps{M})} \ (= \ \mu_M)$ and
$- t  \tau_4  =  - c  b_M/a_{K_{\pi_2}^{-}(\eps{M})} \  (= \ \nu_M)$
hold, because $\lambda_M=-\tau_1 a_{K_{\pi_2}^{+}(\eps{M})}/ (\tau_2  a_{K_{\pi_2}^{-}(\eps{M})})$.
With such a $t$, $\exp(t \, \mathrm{ad}x_{\overline{\varepsilon_{M}}})(x_{\eps{M}} +\lambda_M x_{-\eps{M}})=z_M$.
Hence $z_M$ is semisimple since $x_{\eps{M}}+\lambda_M x_{-\eps{M}}$ is.
\end{proof}

We can now complete the proof of Theorem~\ref{t-lac}:

\medskip

Let $(\underline{a},\underline{b})\in (\C^*)^{\kg{\pi_1}+\kg{\pi_2}}$
such that $\varp{u}{\q_{\pi_1,\pi_2}}$ is $\q_{\pi_1,\pi_2}$-regular where $u=(\underline{a},\underline{b})$.
The orthogonal of $E_{\pi_1,\pi_2}$ in $\mathfrak{h}$ is contained in $\q_{\pi_1,\pi_2}(\var{u})$.
Then, by Lemma \ref{l-yzt}, it suffices to prove that the elements $y_K, z_M, t_N$,
for $K \in \cali{K}_{\pi_1} \cap \cali{K}_{\pi_2}$, $M \in \widetilde{\cali{K}}_1^{(2)}$,
and $N \in \widetilde{\cali{K}}_2^{(1)}$, are linearly independent.
Indeed, if so, the stabilizer of $\varp{u}{\q_{\pi_1,\pi_2}}$ in $\q_{\pi_1,\pi_2}$
contains a (commutative) subalgebra which consists of semisimple elements of $\q_{\pi_1,\pi_2}$ and of dimension
$(\rk\g - \dim E_{\pi_1,\pi_2}) + \#(\cali{K}_{\pi_1} \cap \cali{K}_{\pi_2})+
\# \widetilde{\cali{K}}_1^{(2)}+  \# \widetilde{\cali{K}}_2^{(1)}$.
But the hypothesis of Theorem~\ref{t-lac} tells us that $\#(\cali{K}_{\pi_1} \cap \cali{K}_{\pi_2})+
\# \widetilde{\cali{K}}_1^{(2)}+  \# \widetilde{\cali{K}}_2^{(1)}=\dim (E_{\pi_1} \cap E_{\pi_2})=
\kg{\pi_1}+\kg{\pi_2}-\dim E_{\pi_1,\pi_2}$ (cf.~Definition~\ref{d-lac}).
Hence, by formula~(\ref{eq-ind}), this subalgebra is the stabilizer of $\varp{u}{\q_{\pi_1,\pi_2}}$.

Now, by construction, $M,N,K_{\pi_2}^{+}(\eps{M})$, $K_{\pi_2}^{-}(\eps{M})$, $
K_{\pi_1}^{-}(\eps{N})$, $K_{\pi_1}^{+}(\eps{N})$
do not belong to $\cali{K}_{\pi_1} \cap \cali{K}_{\pi_2}$.
Moreover, $M \in \cali{K}_{\pi_1} \setminus \cali{K}_{\pi_2}$ and $N\in \cali{K}_{\pi_2} \setminus \cali{K}_{\pi_1}$,
whence the expected statement.

\section{Non quasi-reductive parabolic subalgebras}\label{S-4}

So far, our results (Theorem~\ref{t-lac}) only provide examples of quasi-reductive parabolic subalgebras.
It is much trickier to prove that a given Lie algebra is \emph{not} quasi-reductive.
Indeed, to prove that a given parabolic subalgebra is quasi-reductive,
one can make explicit computations, cf.~Section~\ref{S-5}.
In this section we exhibit examples of non quasi-reductive parabolic subalgebras.

\subsection{}

We first discuss the case of the parabolic subalgebras $\p_{\pi'}^{+}$ where $\pi'$ only consists of one simple root.
For $\alpha \in \pi$, denote the parabolic subalgebra $\p_{\{\alpha\}}^{+}$ simply by $\p_{\alpha}^+$\index{$\p_{\alpha}^+$}.
Thanks to Theorem~\ref{t-one} we have a criterion for the quasi-reductivity of $\p_{\alpha}^{+}$:

\begin{theorem} \label{t-one}
Let $\alpha$ be in $\pi$.
Then the parabolic subalgebra $\p_{\alpha}^{+}$ is quasi-reductive if
and only if one of the following two conditions holds:
$\alpha\in\widetilde{\Delta}^{+}_{\pi}$ or $\{\alpha\}\cup \cali{E}_{\pi}$ consists of linearly independent elements.
\end{theorem}

If one of the above two conditions are satisfied, then the cascades of $\{\alpha\}$ and of $\pi$ are well-interlaced;
so, it is clear that $\p_{\alpha}^+$ is quasi-reductive by Theorem~\ref{t-lac}.
Thus, Theorem~\ref{t-one} provides a converse to Theorem~\ref{t-lac} for $\pi_1=\{\alpha\}$ and $\pi_2=\pi$.

\begin{proof}
We only need to show that if $\p_{\alpha}^+$ is quasi-reductive then $\alpha$ satisfies
one of the two conditions of the theorem.
Suppose that $\p_{\alpha}^+$ is quasi-reductive.
If $\alpha$ does not satisfy any of the above conditions, then $\alpha\in E_{\pi}$, and $\alpha$ is not
an element of $\widetilde{\Delta}_{\pi}^{+} \sqcup \cali{E}_{\pi}$.
By Proposition~\ref{p-Zar}, we can find $w$ in $\p_{\alpha}^-$
such that $\varp{w}{\p_{\alpha}^+}$ is regular and of reductive type for $\p_{\alpha}^+$.
Moreover, by Proposition~\ref{p-bor}(ii), since $\alpha \not\in\cali{E}_{\pi}$, we can suppose that $w$ is of the form:
$w=a x_{-\alpha} + h + b x_{\alpha} +\uq{}$ with $a,b \in \C$, $h \in \mathfrak{h}$.
Let us remind that the stabilizer of
$\varp{\uq{}}{\b^+}$ in $\mathfrak{b}^+$  is the orthogonal of $E_{\pi}$ in $\h$ (Proposition \ref{p-bor}(i)).
Consequently, as $\alpha \in E_{\pi}$, we have $[\b^+(\var{\uq{}}),w]=\{0\}$,
whence $\b^+(\var{\uq{}}) \subset \p_{\alpha}^+(\var{w})$.
In addition, by formula~(\ref{eq-ind}), $\ind \p_{\alpha}^+= \ind \b^+ +1$.
So, $\b^+(\var{\uq{}})$ is an  hyperplane of $\p_{\alpha}^+(\var{w})$ (cf.~\cite[Lemma 4.5]{TY2}).
Now choose $x$  in $\p_{\alpha}^+$ such that the decomposition
\begin{eqnarray}\label{eq-one}
\p_{\alpha}^+(\var{w})= \C x \oplus \bigcap\limits_{K \in \cali{K}_{\pi}} \ker\eps{K}
\end{eqnarray}
holds.
By the choice of $w$, $\p_{\alpha}^+(\var{w})$ is an abelian Lie algebra consisting of semisimple elements.
In particular $x$ must be semisimple.
Write the element $x$ as follows: $x=\lambda  x_{-\alpha} +h' + \mu  x_{\alpha}+x^{+}$
with $\lambda, \mu \in \C$, $h' \in \mathfrak{h}$ and $x^{+} \in \m_{\alpha}^+$.
From the fact $[x,w] \in \m_{\alpha}^+$, we deduce that
$h' \in \bigcap\limits_{K \in \cali{K}_{\pi}} \ker\eps{K}$.
So we can assume that $h'=0$ according to (\ref{eq-one}).
Hence $\lambda \mu\not=0$, since $x$ is semisimple.

Since $\alpha$ is not in $\cali{E}_{\pi}$, $\eps{K_{\pi}^{+}(\alpha)}-\alpha$ is a (positive) root.
In turn, suppose that $\eps{K} -\alpha$ is a root, for $K\in\cali{K}_{\pi}$..
As $\alpha$ is a simple root, $\eps{K} -\alpha$ is necessarily a positive root, so $K=K_{\pi}^+(\alpha)$.
Therefore, we have
$$[x,w] = \lambda \sum\limits_{L \in \cali{K}_{\pi}^{-}(\alpha)} [x_{-\alpha},x_{-\eps{L}}]
+ \mu  [x_{\alpha},x_{-\eps{K_{\pi}^{+}(\alpha)}}] + \lambda \alpha(h) x_{-\alpha}
+ (a \mu-  b\lambda) h_{\alpha} - \mu \alpha(h) x_{\alpha} + [x^{+},w]$$
As $[x,w]\in \m_{\alpha}^+$, the bracket $[x_{-\eps{K_{\pi}^{+}(\alpha)}},x_{\alpha}]$ must be compensated.
This bracket cannot be compensated by the term $[x^{+},w]$.
Indeed, if this were the case, then there  would exist $K \in \mathcal{K}_{\pi}$ and $\beta \in\Delta_{\pi}^{+} \setminus \{\alpha\}$
such that $\varepsilon_{K_{\pi}^{+}(\alpha)}-\alpha=\varepsilon_K - \beta$.
But this would force $K=K_{\pi}^{+}(\alpha)$ and so $\alpha=\beta$, which is impossible.
We deduce that there is $L \in \mathcal{K}_{\pi}^{-}(\alpha)$ such that $\varepsilon_{K_{\pi}^{+}(\alpha)}-\alpha=\varepsilon_L+\alpha$.
Thus, $\alpha=\frac{1}{2}(\varepsilon_{K_{\pi}^{+}(\alpha)}-\varepsilon_L)$
that is $\alpha \in \widetilde{\Delta}_{\pi}^{+}$ which contradicts our assumption on $\alpha$.
\end{proof}

According to Theorem~\ref{t-one}, we list the simple roots $\alpha$ corresponding to a
\emph{non} quasi-reductive parabolic subalgebra $\mathfrak{p}_{\alpha}^{+}$
(for simple $\g$) in Table~\ref{T-one}.

{
\begin{table}
\tiny 
\begin{tabular}{|c|c|c|c|c|c|c|}
\hline
&&&&&&\\
B$_{\ell}$, $\ell \geq 3$ & D$_{\ell}$, $\ell \geq 4$ & G$_2$ & F$_4$ & E$_6$ & E$_7$ & E$_8$ \\
&&&&&&\\
\hline
&&&&&&\\
$\alpha_i$, $2 \leq i \leq \ell-1$, $i$ even & $\alpha_i$, $2 \leq i \leq \ell-2$, $i$ even
& $\alpha_1$ & $\alpha_1$ & $\alpha_2$ & $\alpha_1$, $\alpha_4$, $\alpha_6$ &
$\alpha_1$, $\alpha_4$, $\alpha_6$, $\alpha_8$\\
&&&&&&\\
\hline
\end{tabular}
\vspace{.3cm}
\caption{\label{T-one}
The parabolic subalgebras $\p_{\alpha}^+$ which are not quasi-reductive.}
\end{table}}

\begin{remark}
In the exceptional case, Table \ref{T-one} shows that there is always at least one non quasi-reductive parabolic subalgebra.
\end{remark}

\subsection{}

We now exhibit a few more parabolic subalgebras which are not quasi-reductive (Theorem~\ref{t-nE78}
and Theorem~\ref{t-nE6}), all in type E.

\begin{theorem}\label{t-nE78}
{\rm (i)} If $\g$ is of type $\mathrm{E}_7$ and if $\pi'$ is one of the subsets
$\{\alpha_1,\alpha_3,\alpha_4\}$, $\{\alpha_4,\alpha_5,\alpha_6\}$, or $\{\alpha_1,\alpha_3,\alpha_4,\alpha_5,\alpha_6\}$,
then $\p_{\pi'}^{+}$ is not quasi-reductive.

{\rm (ii)} If $\g$ is of type $\mathrm{E}_8$ and if $\pi'$ is one of the subsets
$\{\alpha_1,\alpha_3,\alpha_4\}$,  $\{\alpha_4,\alpha_5,\alpha_6\}$, $\{\alpha_6,\alpha_7,\alpha_8\}$,
$\{\alpha_1,\alpha_3,\alpha_4,\alpha_5,\alpha_6\}$
$\{\alpha_4,\alpha_5,\alpha_6,\alpha_7,\alpha_8\}$,
or $\{\alpha_1,\alpha_3,\alpha_4,\alpha_5,\alpha_6,\alpha_7,\alpha_8\}$,
then $\p_{\pi'}^{+}$ is not quasi-reductive.
\end{theorem}

The indices of the parabolic subalgebras considered in Theorem~\ref{t-nE78} are given in Table~\ref{T-E8}.
Note that for $\g$ of type $\mathrm{E}_7$ or $\mathrm{E}_8$, and $\pi'=\{\alpha_4,\alpha_5,\alpha_6\}$,
$\p_{\pi'}^+$ is not quasi-reductive by Theorem~\ref{t-tran} and Example~\ref{ex-d6}.\\

In the proof of the theorem and in Lemma~\ref{l-nE78} below, we make use of the following notations:
If $\pi'$ is a connected subset of $\pi$, $\widetilde{\pi}'$ is defined to be
the connected subset of $\pi'$ satisfying $\cali{K}_{\pi'}=\{\pi'\} \cup \cali{K}_{\widetilde{\pi}'}$ and $u_{\pi'}^+$ is the element
$\sum\limits_{\eps{}\in\cali{E}_{\pi} \setminus\Delta_{\pi'}^+} x_{\eps{}}$..
Note that the element $u_{\pi'}^+ + u_{\pi'}^-$ is a semisimple element of $\g$.
Assume that $\g$ is of type E$_8$. Set:

\medskip

$\begin{array}{l}
\alpha_{11}=\alpha_3+\alpha_4, \ \alpha_{12}=\alpha_4+\alpha_5, \
\alpha_{13}=\alpha_5+\alpha_6, \ \alpha_{14}=\alpha_6+\alpha_7,\\
\alpha_{19}=\alpha_3+\alpha_4+\alpha_5, \ \alpha_{20}=\alpha_4+\alpha_5+\alpha_6, \ \alpha_{21}=\alpha_5+\alpha_6+\alpha_7,\\
\alpha_{27}=\alpha_3+\alpha_4+\alpha_5+\alpha_6, \ \alpha_{28}=\alpha_4+\alpha_5+\alpha_6+\alpha_7,
\alpha_{35}=\alpha_3+\alpha_4+\alpha_5+\alpha_6+\alpha_7
\end{array}$

\medskip
\noindent
and denote by $I_{\widetilde{\pi}'}$ the set of integers $i$ such that $\alpha_i \in \Delta_{\widetilde{\pi}'}$.
Whenever $\alpha_i$ is defined, $x_i$ and $y_i$ stand for $x_{\alpha_i}$ and $x_{-\alpha_i}$ respectively.
Consider the following equations:

\medskip

\begin{tabular}{llcccllcc}
(E1) & $\mu_{4}+\nu_{19}$ &=&0  &\hspace{2cm}& \hspace{2cm}(G1) & $\mu_{11}-\nu_{12}$ &=&0\\
(F1) & $\mu_{19}+\nu_{4}$ &=&0  &\hspace{2cm}& \hspace{2cm}(H1) & $\mu_{12}-\nu_{11}$ &=&0\\
(E2) & $-\mu_{6}+\nu_{21}$ &=&0 &\hspace{2cm}& \hspace{2cm}(G2) & $\mu_{13}+\nu_{14}$ &=&0\\
(F2) & $\mu_{21}-\nu_{6}$ &=&0  &\hspace{2cm}& \hspace{2cm}(H2) & $\mu_{14}+\nu_{13}$ &=&0\\
(E3) & $\mu_{20}+\nu_{35}$ &=&0 &\hspace{2cm}& \hspace{2cm}(G3) & $\mu_{27}-\nu_{28}$ &=&0\\
(F3) & $\mu_{35}+\nu_{20}$ &=&0 &\hspace{2cm}& \hspace{2cm}(H3) & $\mu_{28}-\nu_{27}$ &=&0\\
\end{tabular}

\medskip
\noindent
in the variables $\mu_i$ and $\nu_i$. 
Set $\pi'_1=\{\alpha_1,\alpha_3,\alpha_4,\alpha_5,\alpha_6\}$,
$\pi'_2=\{\alpha_4,\alpha_5,\alpha_6,\alpha_7,\alpha_8\}$
and $\pi'_3=\{\alpha_1,\alpha_3,\alpha_4,\alpha_5,\alpha_6,\alpha_7,\alpha_8\}$.
We now introduce subspaces $\mathfrak{a}_k$ of $\g_{\widetilde{\pi_k}'}$ (for $k=1,2,3$) as
follows: \\
- for $k=1,2$, $\mathfrak{a}_k$ is the space of elements $\sum\limits_{\eps{} \in \cali{E}_{\pi} \cap \widetilde{\pi}'_k} \lambda_{\eps{}} h_{\eps{}}
+ \sum\limits_{i\in I_{\widetilde{\pi}'_k}} (\mu_i x_i +\nu_i y_i)$ with
$(\lambda_{\eps{}})_{\eps{}\in \cali{E}_{\pi} \cap \widetilde{\pi}'_k}$ in $\C^{|\cali{E}_{\pi} \cap \widetilde{\pi}'_k|}$
and where $((\mu_i)_{i \in I_{\widetilde{\pi}'_k}},(\nu_j)_{j\in I_{\widetilde{\pi}'_k}})$ run through the set of the solutions of
the homogeneous linear system defined by the equations (E$k$), (F$k$), (G$k$), (H$k$).\\
- $\mathfrak{a}_3$ is the space of elements $\sum\limits_{\eps{} \in \cali{E}_{\pi} \cap \widetilde{\pi}'_3} \lambda_{\eps{}} h_{\eps{}}
+ \sum\limits_{i\in I_{\widetilde{\pi}'_3}} (\mu_i x_i +\nu_i y_i)$ with
$(\lambda_{\eps{}})_{\eps{}\in \cali{E}_{\pi} \cap \widetilde{\pi}'_3}$ in $\C^{|\cali{E}_{\pi} \cap \widetilde{\pi}'_3|}$
and where $((\mu_i)_{i\in I_{\widetilde{\pi}'_3}},(\nu_j)_{j\in I_{\widetilde{\pi}'_3}})$ runs through the set of the solutions of
the homogeneous linear system defined by all twelve equations.

Here is a technical lemma used in the proof of Theorem~\ref{t-nE78}:

\begin{lemma} \label{l-nE78}
Assume that $\g$ is of type $\mathrm{E}_8$.
Then $\mathfrak{a}_k$, for $k=1,2,3$, is the centralizer in $\g_{\widetilde{\pi}'_k}$ of the semisimple element $u_{\pi'}^+ + u_{\pi'}^-$.
It is a reductive Lie algebra and its rank is at most $\ind\p_{\pi'}^+ -1$.
\end{lemma}

\begin{proof} Let $k \in\{1,2,3\}$.
The fact that $\mathfrak{a}_k$ centralizes $u_{\pi'}^-$ can be checked without difficulty.
As $\mu$ and $\nu$ play the same role in the equations (E$k$), (F$k$), (G$k$), (H$k$), we deduce that
$\mathfrak{a}_k$ centralizes $u_{\pi'}^+$ too; hence $\mathfrak{a}_k$ centralizes $u_{\pi'}^+ +u_{\pi'}^-$.
Then $\mathfrak{a}_k$ is a reductive Lie algebra as an intersection between a reductive
Lie algebra and the centralizer in $\g$ of a semisimple element of $\g$.

Next we show: $\rk\mathfrak{a}_k \leq \ind\p_{\pi'}^+ -1$.
We can readily verify from the equations defining $\a_k$ that the center of $\mathfrak{a}_k$ is zero.
Therefore, the rank of $\mathfrak{a}_k$ is strictly smaller that the one of $\g_{\widetilde{\pi}'_k}$.
Indeed, if not, $\mathfrak{a}_k$ is a Levi subalgebra of $\g_{\widetilde{\pi}'_k}$ since $\g_{\widetilde{\pi}'_k}$ has type A.
But any proper Levi subalgebra of $\g_{\widetilde{\pi}'_k}$ has a non trivial center.
So, for $k=1,2$, we get $\rk \mathfrak{a}_k \leq 2$ since $\rk\g_{\widetilde{\pi}'_k}=\ind \p_{\pi'_k}^+=3$
whence the statement.

For $k=3$, what foregoes yields $\rk \mathfrak{a}_3 \leq 4$ since the rank of $\g_{\widetilde{\pi}'_3}$ is 5.
We have to show: $\rk \mathfrak{a}_3 < \ind\p_{\pi'_3}^+=4$.
The space $\mathfrak{a}_3$ has dimension 21.
But there is no reductive Lie subalgebra of rank 4 and of dimension 21 since $21-4$ is not even.
As a result, we get $\rk \mathfrak{a}_3 < 4$.
\end{proof}

Here is the proof of Theorem~\ref{t-nE78}:

\begin{proof}[Proof of Theorem~\ref{t-nE78}]
By the transitivity property (Theorem~\ref{t-tran}), statement (ii) implies (i).
So we only consider the case of $\mathrm{E}_8$.
Let $\pi'$ be one of the subsets as described in (ii).
Assume that $\p_{\pi'}^+$ is quasi-reductive.
We will show that this leads to a contradiction.
Choose $w \in \p_{\pi'}^{-}$ such that the following two conditions are satisfied:

- $\varp{w}{\p_{\pi'}^+}$ is $\p_{\pi'}^+$-regular and of reductive type for $\p_{\pi'}^+$;

- $\varp{w}{\m_{\pi'}^+}$ belongs to the $B$-orbit of $\varp{u_{\pi'}^{-}}{\m_{\pi'}^+}$.\\
This choice of $w$ is possible by Proposition~\ref{p-Zar} and Proposition~\ref{p-bor}(ii).
By the second condition, we can assume that $w=w'+ \uq{\pi'}$ with $w'\in \l_{\pi'}$.
Let $x$ be an element of the stabilizer $\varp{w}{\p_{\pi'}^+}$ in $\p_{\pi'}^+$;
we write $x=h+x'+x^+$, with $h\in\h$, $x'\in \n_{\pi'}^-\oplus \n_{\pi'}^+$ and $x^+ \in \m_{\pi'}^+$.
The fact $[x,w]\in\m_{\pi'}^+$ forces $h \in \bigcap\limits_{\eps{} \in \cali{E}_{\pi} \setminus \Delta_{\pi'}^+}\ker \eps{}$.
From that, we deduce that $h$ belongs to the subspace of $\h$ generated by the elements $h_{\eps{}}$,
for  $\eps{} \in \cali{E}_{\pi} \cap \pi' \subset \widetilde{\pi}'$ (use Table~\ref{T-exc}).
Now for $\alpha \in \Gamma_{\pi'}$, one obtains that $\eps{K_{\pi}^{+}(\alpha)} \not \in \Delta_{\pi'}^+$
and we claim that $x'$ has zero coefficient in $\g_{\alpha}$.
Otherwise, there must be $\beta \in \Delta_{\pi'}^{+}$ and $K \in \cali{K}_{\pi}$ such that
$\alpha - \eps{K_{\pi}^{+}(\alpha)}=- (\beta+ \eps{K})$.
One can check that for each of the subsets $\pi'$ such an equality is not possible (use Table~\ref{T-exc}).
To summarize, we obtain the inclusion:
\begin{eqnarray} \label{eq-dec}
\p_{\pi'}^+(\var{w}) &\subset  &\g_{\widetilde{\pi}'} \oplus \mathfrak{H}_{\pi'}^{-} \oplus\m_{\pi'}^+,
\end{eqnarray}
where $\mathfrak{H}_{\pi'}^{-}$ is the Heisenberg Lie algebra generated by the
$\g_{-\alpha}$, $\alpha \in \Gamma_{\pi'}$.
Let $\mathfrak{t}$ be the image of $\p_{\pi'}^+(\var{w})$ by the projection map
from $\g_{\widetilde{\pi}'} \oplus \mathfrak{H}_{\pi'}^{-} \oplus \m_{\pi'}^+$ to $\g_{\widetilde{\pi}'}$.
As $\p_{\pi'}^+(\var{w})$ is a torus of $\g$ by hypothesis, (\ref{eq-dec}) shows that $\mathfrak{t}$
is a torus of $\g_{\widetilde{\pi}'}$ of dimension $\ind \p_{\pi'}^+=\dim \p_{\pi'}^+(\var{w})$.

For the first three subsets, with $\pi'$ of rank 3, $\p_{\pi'}^+$ has index $2$ but $\g_{\widetilde{\pi}'}$ has rank $1$.
So we get a contradiction.

The remaining cases, with $\pi'$ of rank 5 or 7, require more work.
Let us describe the torus $\t$.
To do that, we consider on one hand the roots $\alpha \in \Delta_{\widetilde{\pi}'}^+$ with $\eps{K_{\pi}^{+}(\alpha)} \not \in \Delta_{\pi'}^+$
for which there exist $\beta \in \Delta_{\pi'}^{+}$ and $K \in \cali{K}_{\pi}$ such that $\alpha - \eps{K_{\pi}^{+}(\alpha)}=- (\beta+ \eps{K})$.
On the other hand, we consider the roots $\alpha \in \Delta_{\widetilde{\pi}'}^-$ for which there is $\eps{} \in \cali{E}_{\pi}\setminus \Delta_{\pi'}$ such that $\alpha +\eps{}$ is a root.
All the possible roots give rise to equations describing $\mathfrak{t}$.
Let $k\in\{1,2,3\}$ and use the notations introduced before Lemma~\ref{l-nE78}.
The equations what we obtained are precisely the equations (E$k$), (F$k$), (G$k$), (H$k$) for $k=1,2$, and
all twelve equations above for $k=3$. 
Thereby $\mathfrak{t}$ is contained in the reductive Lie algebra $\mathfrak{a}_k$.
But the torus $\t$ has dimension $\ind \p_{\pi'}^+$ and this contradicts Lemma~\ref{l-nE78}.
\end{proof}

\begin{remark} \label{r-tor}
Proceeding with the proof of Lemma~\ref{l-nE78}, one readily obtains that $\a_k$, for $k=1,2,3$, has precisely dimension $\ind \p_{\pi'}^+ - 1$
(note that $\dim \a_1=\dim \a_2 = 10$ and $\dim \a_3 = 21$).
Then, the proof of Theorem~\ref{t-nE78} shows that the dimension of the torus part of generic stabilizers is $\ind \p_{\pi'}^+ - 1$.
This dimension is given, for each case, in the last column of Table~\ref{T-E8}.
\end{remark}

We end the section with an example of non quasi-reductive parabolic subalgebra in E$_6$.
As noticed in Remark~\ref{r-add}, Theorem~\ref{t-nE6} shows that the additivity property fails in type E$_6$:

\begin{theorem} \label{t-nE6}
If $\g$ if of type $\mathrm{E}_6$ and if $\pi'=\{\alpha_1,\alpha_2,\alpha_3,\alpha_4,\alpha_6\}$,
then $\p_{\pi'}^{+}$ is not quasi-reductive.
\end{theorem}

By symmetry, if $\pi'=\{\alpha_1,\alpha_2,\alpha_4,\alpha_5,\alpha_6\}$, then $\p_{\pi'}^{+}$ is not quasi-reductive, either.

\begin{proof}
Choose $w \in \p_{\pi'}^{-}$ such that the following two conditions are satisfied:

- $\varp{w}{\p_{\pi'}^+}$ is $\p_{\pi'}^+$-regular;

- $\varp{w}{\n^+}$ belongs to the $B$-orbit of $\varp{u^{-}}{\n^+}$.\\
This choice is possible by Proposition~\ref{p-Zar} and Proposition~\ref{p-bor}(ii).
By the second condition, we can assume that $w=w'+ \uq{}$ where $w'$ is in $\h \oplus \n_{\pi'}^+$.
For $x \in \p_{\pi'}^+$, we write $x=h+x'+x^+$ with $h\in\h$, $x'\in \n_{\pi'}^-\oplus \n_{\pi'}^+$ and $x^+ \in \m_{\pi'}^+$.
Set:

\medskip

$\begin{array}{l}
\alpha_{7}=\alpha_1+\alpha_3, \ \alpha_{8}=\alpha_2+\alpha_4, \ \alpha_{9}=\alpha_3+\alpha_4,\\
\alpha_{12}=\alpha_1+\alpha_3+\alpha_4, \ \alpha_{13}=\alpha_2+\alpha_3+\alpha_4, \\
\alpha_{17}=\alpha_1+\alpha_2+\alpha_3+\alpha_4
\end{array}$

\medskip

\noindent
and let $I_{\pi'}$ be the set of integers $i$ such that $\alpha_i\in \Delta_{\pi'}$.
Then, for $i\in I_{\pi'}$, $x_i$, $y_i$ and $h_i$ stand for $x_{\alpha_i}$, $x_{-\alpha_i}$ and $h_{\alpha_i}$ respectively.
Write $x'=\sum\limits_{i\in I_{\pi'}} \mu_i x_i + \sum\limits_{i=1}^{6} \lambda_i h_i
+\sum\limits_{i\in I_{\pi'}} \nu_i y_i$ and $w'=h_0+\sum\limits_{l\in I_{\pi'}} a_l x_l$ with $h_0 \in\h$
and $(\mu_i,\lambda_j,\nu_k,a_l)_{i,j,k,l} \in \C^{3|I_{\pi'}|+6}$.

From $[x,w]\in\m_{\pi'}^+$, we first deduce that
$h$ belongs to $\ker \eps{}$ for any $\eps{} \in \cali{E}_{\pi} \setminus \Delta_{\pi'}^+$
whence we get $\lambda_1=-\lambda_6$ and $\lambda_3=-\lambda_5$.
Next, we argue as at the end of Theorem~\ref{t-nE78}(ii):
we use the roots $\alpha \in \Delta_{\pi'}^+$ such that $\eps{K_{\pi}^{+}(\alpha)} \not \in \Delta_{\pi'}^+$
and for which there exist $\beta \in \Delta_{\pi'}^{+}$ and $K \in \cali{K}_{\pi}$ such that
$\alpha - \eps{K_{\pi}^{+}(\alpha)}=- (\beta+ \eps{K})$.
This enables us to show that $\mu_i=0$ for any $i\in I_{\pi'} \setminus\{1,4,6\}$ and that $\nu_6=\mu_1$, $\mu_6=\nu_1$.
Now, we consider the terms in $x_{\alpha}$
for $\alpha\in\Delta_{\pi'}$ and in $h_{\alpha}$ for $\alpha\in \pi'$ of $[x,w]$.
All these terms have to be zero; this gives us equations.
Some of them involve the terms in $x_{\alpha}$ for certain $\alpha \in \Delta_{\pi}^+\setminus\Delta_{\pi'}^+$ but we can
eliminate these variables and obtain equations whose variables are only
the $(\lambda_i)_i$'s, $(\mu_j)_j$'s, and $(\nu_k)_k$'s, for $i=1,3,4$, $j=1,4$  and $k\in I_{\pi'}$.
Here are these equations:

\medskip

\begin{tabular}{llcc}
(X1)& $2a_1 \lambda_1 - a_1 \lambda_3 +(\alpha_6-\alpha_6)(h_0) \mu_1+a_7\nu_3+a_{12}\nu_{9} + a_{17}\nu_{13}$& =&0\\
(X2)& $-a_2\lambda_4 + a_8 \nu_4 +2\nu_{8}+a_{13}\nu_{9}+a_{17}\nu_{12}$&=&0\\
(X3)& $-a_3\lambda_1+2a_3\lambda_3-a_3\lambda_4 -a_7\nu_1 +a_9\nu_4 +a_{13}\nu_{8}$&=&0\\
(X4)& $2a_4\lambda_4-\alpha_4(h_0)\mu_4-a_8\nu_2-a_9\nu_3-a_{12}\nu_7$&=&0\\
(X5)& $-2\lambda_4 + \alpha_4(h_0)\nu_4 -a_2\nu_{8}-a_3\nu_{9}-a_7\nu_{12}$&=&0\\
(X6)& $-2a_6\lambda_1+a_6\lambda_3+(\alpha_1-\alpha_6)(h_0)\nu_1+a_3\nu_7+a_9\nu_{12}+a_{13}\nu_{17}$&=&0\\
(X7)& $a_7\lambda_1+a_7\lambda_3-a_7\lambda_4-a_3\mu_1+a_{12}\nu_4+a_{17}\nu_{8}$&=&0\\
(X8)& $a_8\lambda_4+a_2\mu_4-2\nu_2-a_{13}\nu_3-a_{17}\nu_7$&=&0\\
(X9)& $-a_9\lambda_1+2a_9\lambda_3+a_9\lambda_4+a_3\mu_4-a_{12}\nu_1-a_{13}\nu_2$&=&0\\
(X12)& $a_{12}\lambda_1+a_{12}\lambda_3+a_{12}\lambda_4-a_9\mu_1+a_7\mu_4-a_{17}\nu_2$&=&0\\
(X13)& $-a_{13}\lambda_1-a_{17}\nu_1$&=&0\\
(X17)& $a_{17}\lambda_1+a_{17}\lambda_3-a_{13}\mu_1$&=&0\\
(H1)&  $a_6\mu_1-a_1 \nu_1-a_7\nu_7-a_{12}\nu_{12}-a_{17}\nu_{17}$&=&0\\
(H3)& $-a_3 \nu_3-a_7\nu_7-a_9\nu_{9}-a_{12}\nu_{12}-a_{13}\nu_{13}-a_{17}\nu_{17}$&=&0\\
(H4)& $-2\mu_4-a_2\nu_2+2a_4\nu_4+a_8\nu_{8}+2a_9\nu_{9}+2a_{12}\nu_{12}+a_{13}\nu_{13}+a_{17}\nu_{17}$&=&0\\
\end{tabular}

\medskip

Using a computer algebra system, we show that
for any $((a_i)_{i \in I_{\pi'}},(\alpha_i(h_0))_{i\in\pi'})$ in an open dense subset of $\C^{|I_{\pi'}|}\times \C^6$,
the above homogeneous linear system has rank 14, $a_{13} a_{17}\not= 0$, and any of its solution
$((\lambda_i)_{i=1,3,4},(\mu_j)_{j=1,4},(\nu_k)_{k\in I_{\pi'}})$ verifies $\lambda_3=0$.
We can (and do) assume that $((a_i)_{i \in I_{\pi'}},(\alpha_i(h_0))_{i\in\pi'})$ belongs to this open subset;
in particular $a_{13} a_{17}\not= 0$.
From the equations (X13) and (X17), we obtain that any solution of this system verifies $\lambda_1^2+ \mu_1 \nu_1=0$
because $\lambda_3=0$.
Since $\mu_i=0$ for any $i\in I_{\pi'}\setminus\{1,4,6\}$ as observed previously,
this shows that $x'$ is a nilpotent element of $\l_{\pi'}$; so $x$ is a nilpotent element of $\g$.
As a consequence, $\p_{\pi'}^+$ is not quasi-reductive.
\end{proof}

%
\section{Explicit computations and classification} \label{S-5}
%

We assume in this part that $\g$ is simple of exceptional type.
Together with Theorem~\ref{t-dkt}, the next two theorems (Theorem~\ref{t-EF} and
Theorem~\ref{t-E6}) complete the classification of quasi-reductive parabolic subalgebras of simple Lie algebras.
The goal of this section is to prove these theorems.

\begin{theorem}\label{t-EF}
Assume that $\g$ is of type $\mathrm{G}_2$, $\mathrm{F}_4$, $\mathrm{E}_7$ or $\mathrm{E}_8$.
Let $\pi'$ be a subset of $\pi$.

{\rm (i)} If $\g$ is of type $\mathrm{G}_2$, then $\p_{\pi'}^+$ is quasi-reductive
if and only if $\pi'$ is different from $\{\alpha_1\}$.

{\rm (ii)} If $\g$ is of type $\mathrm{F}_4$, then $\p_{\pi'}^+$ is quasi-reductive
if and only if each connected component of $\pi'$ is different from $\{\alpha_1\}$.

{\rm (iii)} If $\g$ is of type $\mathrm{E}_7$, then $\p_{\pi'}^+$ is quasi-reductive
if and only if each connected component of $\pi'$ is different from the subsets
$\{\alpha_1\}$, $\{\alpha_4\}$, $\{\alpha_6\}$, $\{\alpha_1,\alpha_3,\alpha_4\}$,
$\{\alpha_4,\alpha_5,\alpha_6\}$ and $\{\alpha_1,\alpha_3,\alpha_4,\alpha_5,\alpha_6\}$.

{\rm (iv)} If $\g$ is of type $\mathrm{E}_8$, then $\p_{\pi'}^+$ is quasi-reductive
if and only if each connected component of $\pi'$ is different from the subsets
$\{\alpha_1\}$, $\{\alpha_4\}$, $\{\alpha_6\}$, $\{\alpha_8\}$,
$\{\alpha_1,\alpha_3,\alpha_4\}$, $\{\alpha_4,\alpha_5,\alpha_6\}$, $\{\alpha_6,\alpha_7,\alpha_8\}$,
$\{\alpha_1,\alpha_3,\alpha_4,\alpha_5,\alpha_6\}$, $\{\alpha_4,\alpha_5,\alpha_6,\alpha_7,\alpha_8\}$
and $\{\alpha_1,\alpha_3,\alpha_4,\alpha_5,\alpha_6,\alpha_7,\alpha_8\}$.
\end{theorem}

\begin{theorem}\label{t-E6}
Assume that $\g$ is of type $\mathrm{E}_6$ and let $\pi'$ be a subset of $\pi$.
Then $\p_{\pi'}^+$ is quasi-reductive except in the following three cases:

{\rm 1)} $\{\alpha_2\}$ is a connected component of $\pi'$;

{\rm 2)} $\pi'=\{\alpha_1,\alpha_2,\alpha_3,\alpha_4,\alpha_6\}$;

{\rm 3)} $\pi'=\{\alpha_1,\alpha_2,\alpha_4,\alpha_5,\alpha_6\}$.
\end{theorem}

Table~\ref{T-E8} and Table~\ref{T-E6} below summarize the results of
Theorems~\ref{t-EF} and~\ref{t-E6} ;
indeed, whenever $\rk\g=\kg{\pi}$, only the cases where $\pi'$ is connected need to be dealt with thanks to Theorem~\ref{t-add}.
In these tables, the last column gives the dimension of the torus part of a generic stabilizer;
we refer to Remark~\ref{r-tor} for explanations in the types E$_7$ and E$_8$.
For the type E$_6$, let us roughly explain our computations :
in most cases, the subspaces
$\bigcap\limits_{\eps{} \in \cali{E}_{\pi} \cup \Delta_{\pi'}} \ker \eps{}$ of $\h$
yield elements of the generic stabilizers of the regular linear forms of the form
$\varp{w'+u^-}{\p_{\pi'}^+}$ with $w'\in\l_{\pi'}$.
For the cases $\{\alpha_1,\alpha_2,\alpha_6\}$,  $\{\alpha_2,\alpha_3,\alpha_5\}$ and
$\{\alpha_1,\alpha_2,\alpha_3,\alpha_5,\alpha_6\}$,
one can show that the generic stabilizers of these forms also contain
nonzero semisimple elements which do not belong
to $\h$.
Since this is not a central point for our work, we omit the details.

{\begin{table}[h] \tiny
\begin{tabular}{l}

Type F$_4$: \begin{tabular}{|l|c|c|}
\hline
$\pi'$ & $\ind\p_{\pi'}^+$ & dim.~of torus part\\
\hline
$\{\alpha_1\}$ & 1 & 0\\
\hline
\end{tabular}\\
\\
Type E$_7$: \begin{tabular}{|l|c|c|}
\hline
$\pi'$ & $\ind\p_{\pi'}^+$ & dim.~of torus part\\
\hline
$\{\alpha_1\}$ & 1 & 0\\
$\{\alpha_4\}$ &  1 & 0\\
$\{\alpha_6\}$ & 1 & 0 \\
$\{\alpha_1,\alpha_3,\alpha_4\}$ & 2 & 1\\
$\{\alpha_4,\alpha_5,\alpha_6\}$ & 2 & 1\\
$\{\alpha_1,\alpha_3,\alpha_4,\alpha_5,\alpha_6\}$ & 3 & 2\\
\hline
\end{tabular}\\
\\
Type E$_8$: \begin{tabular}{|l|c|c|}
\hline
$\pi'$ & $\ind\p_{\pi'}^+$ & dim.~of torus part\\
\hline
$\{\alpha_1\}$ & 1 & 0\\
$\{\alpha_4\}$ &  1 & 0\\
$\{\alpha_6\}$ & 1 & 0 \\
$\{\alpha_8\}$ & 1 & 0\\
$\{\alpha_1,\alpha_3,\alpha_4\}$ & 2 & 1\\
$\{\alpha_4,\alpha_5,\alpha_6\}$ & 2 & 1\\
$\{\alpha_6,\alpha_7,\alpha_8\}$ & 2 & 1\\
$\{\alpha_1,\alpha_3,\alpha_4,\alpha_5,\alpha_6\}$ & 3 & 2\\
$\{\alpha_4,\alpha_5,\alpha_6,\alpha_7,\alpha_8\}$ & 3 & 2\\
$\{\alpha_1,\alpha_3,\alpha_4,\alpha_5,\alpha_6,\alpha_7,\alpha_8\}$ & 4 & 3\\
\hline
\end{tabular}
\end{tabular}
\vspace{.3cm}
\caption{\label{T-E8}
The non quasi-reductive parabolic subalgebras $\p_{\pi'}^+$
with connected $\pi'$ in F$_4$, E$_7$ and E$_8$  and their indices.}
\end{table}}

{\begin{table}[h] \tiny
\begin{tabular}{|l|c|c|}
\hline
$\pi' \subset \pi$, \ $\pi$ of type E$_6$ & $\ind\p_{\pi'}^+$ & dim.~of torus part\\
\hline
$\{\alpha_2\}$ & 3 & 2\\
$\{\alpha_1,\alpha_2\}$, $\{\alpha_2,\alpha_6\}$ & 2 & 1\\
$\{\alpha_3,\alpha_2\}$, $\{\alpha_2,\alpha_5\}$ & 2 & 1\\
$\{\alpha_1,\alpha_2,\alpha_5\}$, $\{\alpha_2,\alpha_3,\alpha_6\}$ & 1 & 0\\
$\{\alpha_1,\alpha_2,\alpha_6\}$ & 3 & 2\\
$\{\alpha_2,\alpha_3,\alpha_5\}$ & 3 & 2\\
$\{\alpha_1,\alpha_2,\alpha_3\}$, $\{\alpha_2,\alpha_5,\alpha_6\}$ & 2 & 1\\
$\{\alpha_1,\alpha_2,\alpha_3,\alpha_5\}$, $\{\alpha_2,\alpha_3,\alpha_5,\alpha_6\}$ & 1 & 0\\
$\{\alpha_1,\alpha_2,\alpha_3,\alpha_6\}$, $\{\alpha_2,\alpha_1,\alpha_5,\alpha_6\}$ & 1 & 0\\
$\{\alpha_1,\alpha_2,\alpha_3,\alpha_5,\alpha_6\}$ & 3 & 2\\
$\{\alpha_1,\alpha_2,\alpha_3,\alpha_4,\alpha_6\}$, $\{\alpha_1,\alpha_2,\alpha_4,\alpha_5,\alpha_6\}$ & 1 & 0\\
\hline
\end{tabular}
\vspace{.3cm}
\caption{\label{T-E6}
The non quasi-reductive parabolic subalgebras $\p_{\pi'}^+$ in E$_6$ and their indices.}
\end{table}}

\begin{remark} \label{r-addE6} Theorems~\ref{t-EF} and~\ref{t-E6} confirm what was announced
in Remark~\ref{r-add}:
The only cases where the additivity property fails is for
$\g=\mathrm{E}_6$ and $\pi'=\{\alpha_1,\alpha_2,\alpha_3,\alpha_4,\alpha_6\}$
(where $\{\alpha_1,\alpha_2,\alpha_3,\alpha_4\}$ is not connected to $\{\alpha_6\}$),
resp.~for $\g=\mathrm{E}_6$ and $\pi'=\{\alpha_1,\alpha_2,\alpha_4,\alpha_5,\alpha_6\}$.
\end{remark}

By Theorem~\ref{t-one}, Theorem~\ref{t-nE78} and Theorem~\ref{t-nE6}, in order to prove
Theorems~\ref{t-EF} and~\ref{t-E6}, it is enough to show that if $\pi'$ is different from
the subsets listed in Theorems~\ref{t-EF} and~\ref{t-E6} then $\p_{\pi'}^+$ is quasi-reductive.
This is our goal until the end of the paper.
Recall that $\alpha_{\pi}$ is the simple root connected to the lowest root
in the extended Dynkin diagram.
By Theorem~\ref{t-tran}, we can assume that $\pi'$ contains $\alpha_{\pi}$.
Moreover, whenever $\rk\g=\kg{\pi}$, we can assume that $\pi'$ is connected by Theorem~\ref{t-add}.

The case where $\pi'$ has rank 1 was dealt with in Theorem~\ref{t-one}.
In the next subsection, we study the case where $\pi'$ is connected and of rank $2$.
Then we discus the remaining cases in Subsection~\ref{Ss-GAP}.

\subsection{} \label{Ss-2}
Assume that $\g$ is of type F$_4$, E$_6$, E$_7$ or E$_8$ and
let $\pi'$ be a connected subset of $\pi$ of rank $2$ which contains $\alpha_{\pi}$.
Write $\pi'=\{\alpha_{i_1},\alpha_{i_2}\}$ with $\alpha_{i_2}=\alpha_{\pi}$.
Lemma~\ref{l-2} shows that the roots of $\pi'$ have common properties:

\begin{lemma} \label{l-2}
The subset $\pi'$ has type $\mathrm{A}_2$ and there are four integers $j_0,j_1,j_2,j_3$ in $\{1,\ldots,\kg{\pi}\}$ and
a quadruple $(c_0,c_1,c_2,c_3) \in \C^4$ such that the following  properties are satisfied:\\
$\alpha_{i_1} =\eps{j_1}$, $\alpha_{i_2}=\frac{1}{2}(\eps{j_0} -\eps{j_1}-\eps{j_2}-\eps{j_3} )$
and $h_{\eps{\pi'}}=\sum_{k=0}^3 c_k h_{\eps{j_k}}$.
\end{lemma}

\begin{proof}
We verify the properties for each type:

\medskip

\noindent
Type F$_4$: $\pi'=\{\alpha_2,\alpha_1\}$,
with $\alpha_1 (= \alpha_{\pi})=\frac{1}{2}(\eps{1}-\eps{4}-\eps{2}-\eps{3})$ and $\alpha_2=\eps{4}$.
Moreover $h_{\alpha_1+\alpha_2}=\frac{1}{2}(h_{\eps{1}}+h_{\eps{4}}-h_{\eps{2}}-h_{\eps{3}})$.\\
Type E$_6$: $\pi'=\{\alpha_4,\alpha_2\}$,
with $\alpha_2(= \alpha_{\pi})=\frac{1}{2}(\eps{1}-\eps{4}-\eps{2}-\eps{3})$ and $\alpha_4=\eps{4}$.
Moreover $h_{\alpha_2+\alpha_4}=\frac{1}{2}(h_{\eps{1}}+h_{\eps{4}}-h_{\eps{2}}-h_{\eps{3}})$.\\
Type E$_7$: $\pi'=\{\alpha_1,\alpha_3\}$,
with $\alpha_1(= \alpha_{\pi})=\frac{1}{2}(\eps{1}-\eps{6}-\eps{2}-\eps{3})$ and $\alpha_3=\eps{6}$.
Moreover $h_{\alpha_1+\alpha_3}=\frac{1}{2}(h_{\eps{1}}+h_{\eps{6}}-h_{\eps{2}}-h_{\eps{3}})$.\\
Type E$_8$: $\pi'=\{\alpha_7,\alpha_8\}$,
with $\alpha_8(= \alpha_{\pi})=\frac{1}{2}(\eps{1}-\eps{5}-\eps{2}-\eps{3})$ and $\alpha_7=\eps{5}$.
Moreover $h_{\alpha_7+\alpha_8}=\frac{1}{2}(h_{\eps{1}}+h_{\eps{5}}-h_{\eps{2}}-h_{\eps{3}})$.
\end{proof}

Recall that there exist $\underline{a}=(a_1,\ldots,a_{\mathbf{k}_{\mathfrak{g}}}) \in (\C^*)^{\kg{\pi}}$ and $b \in \C^*$,
such that the linear form $\varp{u(\underline{a},b)}{\p_{\pi'}^{+}}$ is $\mathfrak{p}_{\pi'}^{+}$-regular.
Since $\eps{\pi'}=\frac{1}{2}(\eps{j_0}+\eps{j_1}-\eps{j_2}-\eps{j_3})$, the element $\eps{j_0}-\eps{\pi'}$ is a positive root.
Denote by $\beta_2$ and $\beta_3$ the two positive roots $\beta_2=(\eps{j_0}-\eps{\pi'})-\eps{j_2}$ and $\beta_3=(\eps{j_0}-\eps{\pi'})-\eps{j_2}$..
For $\underline{\lambda}=(\lambda_2,\lambda_3,\mu_{0},\mu_{1},\mu_{2},\mu_{3},\nu) \in (\C^*)^{7}$,
we set $x(\underline{\lambda})=x_{-\eps{\pi'}} +  \lambda_2 x_{\beta_2}+\lambda_3 x_{\beta_3} +\sum_{k=0}^3 \mu_{k}  x_{\eps{j_k}} + \nu x_{-\eps{j_1}}$.

\begin{lemma}\label{l2-2}
Let $(\underline{a},b)$ be in $(\C^*)^{\kg{\pi}}\times \C^*$
such that $\varp{u(\underline{a},b)}{\p_{\pi'}^{+}}$ is $\mathfrak{p}_{\pi'}^{+}$-regular.

For a suitable choice of $\underline{\lambda}=(\lambda_2,\lambda_3,\mu_{0},\mu_{1},\mu_{2},\mu_{3},\nu) \in (\C^*)^{7}$,
the element $x(\underline{\lambda})$ lies in the stabilizer of $\varp{u(\underline{a},b)}{\p_{\pi'}^{+}}$ in $\p_{\pi'}^{+}$.
Moreover, for such a $\underline{\lambda}$, we have
$\p_{\pi'}^{+}(\var{u(\underline{a},b)}) = \bigcap\limits_{K \in\cali{K}_{\pi}} \ker \eps{K}\oplus \C x(\underline{\lambda})$
and the element  $x(\underline{\lambda})$ is semisimple.
In particular $\var{u(\underline{a},b)}$ is of reductive type for $\p_{\pi'}^{+}$.
\end{lemma}

\begin{proof}
By definition, we have $\eps{j_2} +\eps{\pi'}= \eps{j_0}-\beta_2$, $\eps{j_3} +\eps{\pi'} = \eps{j_0}-\beta_3$,
$\beta_2 - \eps{j_3} = \beta_3 -\eps{j_2} = \alpha_{i_2}$ and $\eps{\pi'}-\eps{j_1}=\alpha_{i_2}$.
We define the structure constants $\tau_1,\tau_2,\tau_3,\tau_4,\tau_5,\tau_6,\tau_0$
by the following equations:

\begin{center}
\begin{tabular}{ll}
$[ x_{\beta_2},x_{-\eps{j_0}}] = \tau_1  x_{-(\eps{\pi'}+\eps{j_2})}$ \ ; &
$[ x_{-\eps{\pi'}} ,x_{-\eps{j_2}} ] = \tau_2  x_{-(\eps{\pi'}+\eps{j_2})}$\\
$[x_{\beta_3},x_{-\eps{j_0} } ]= \tau_3  x_{-(\eps{\pi'}+\eps{j_3})}$ \ ; &
$[x_{-\eps{\pi'}},x_{-\eps{j_3}} ] = \tau_4  x_{-(\eps{\pi'}+\eps{j_3})}$\\
$[ x_{\beta_2},x_{-\eps{j_3}} ] = \tau_5  x_{\alpha_{i_2}}$  \ ; &
$[ x_{\beta_3},x_{-\eps{j_2}} ] = \tau_6  x_{\alpha_{i_2}}$ \ ;\\
$[x_{-\eps{j_1}},x_{\eps{\pi'}}] = \tau_0 x_{\alpha_{i_2}}$. &
\end{tabular}
\end{center}

\noindent
Set $u=u(\underline{a},b)$ and $x=x(\underline{\lambda})$. We have:
\begin{eqnarray*}
[x,u] =  b
  [x_{-\eps{\pi'}},x_{\eps{\pi'}}] + \sum_{k=0}^3  \mu_{k}  a_{j_k} [x_{\eps{j_k}}, x_{-\eps{j_k}}]
  + a_{j_2}  [x_{-\eps{\pi'}},x_{-\eps{j_2}}] + a_{j_3}  [x_{-\eps{\pi'}},x_{-\eps{j_3}}]
  + \ \lambda_3 a_{j_0} [x_{\beta_3},x_{-\eps{j_0}}] \\+  \lambda_2   a_{j_0} [x_{\beta_2},x_{-\eps{j_0}}]
  + \nu  b [x_{-\eps{j_1}},x_{\eps{\pi'}}] + \lambda_2   a_{j_3} [x_{\beta_2},x_{-\eps{j_3}}] + \lambda_3  a_{j_2}
  [x_{\beta_3},x_{-\eps{j_2}}] +  v
\end{eqnarray*}
where $v$ is in $\mathfrak{m}_{\pi'}^{+}$.
In the above notations, this gives:
\begin{eqnarray*}
[x,u] & = & b (-h_{\eps{\pi'}})+ \sum _{k=0}^3\mu_{k}  a_{j_k} h_{\eps{j_k}}
+  (a_{j_2}  \tau_2+ \lambda_2   a_{j_0}   \tau_1 )  x_{-(\eps{\pi'}+\eps{j_2})}\\
&&+\ (a_{j_3}  \tau_4 + \lambda_3   a_{j_0}  \tau_3 )  x_{-(\eps{\pi'}+\eps{j_3})}
+ (\nu  b  \tau_0+ \lambda_2 a_{j_3} \tau_5 + \lambda_3    a_{j_2}  \tau_6) x_{\alpha_{i_2}} +  v
\end{eqnarray*}
Set $\mu_{k}=(b   c_{k})/a_{j_k}$, for $k=0,1,2,3$.
By Lemma~\ref{l-2} we get $b (-h_{\eps{\pi'}})+ \sum _{k=0,1,2,3} \mu_{k}  a_{j_k} h_{\eps{j_k}}=0$.
Next, we set $\lambda_2 =-a_{j_2}  \tau_2/(a_{j_0}  \tau_1)$ and $\lambda_3 = - a_{j_3}  \tau_4/(a_{j_0}  \tau_3)$
so that the terms in $x_{-(\eps{\pi'}+\eps{j_2})}$ and $x_{-(\eps{\pi'}+\eps{j_3})}$ in $[x,u]$ are both equal to zero.
At last, we choose $\nu$ so that the term in $x_{\alpha_{i_2}}$ in  $[x,u]$ is $0$.
Then the element $x$ stabilizes $\varp{u}{\p_{\pi'}^{+}}$.

Let $\underline{\lambda}$ be as above.
We have thus obtained the inclusion $\bigcap\limits_{K \in\cali{K}_{\pi}} \ker \eps{K} \oplus \C x \subset \p_{\pi'}^{+}(\var{u}) $.
By equation~(\ref{eq-ind}), $\ind \p_{\pi'}^{+}=\rk\g-\kg{\pi}+1$ whence the equality
$\bigcap\limits_{K \in\cali{K}_{\pi}} \ker \eps{K} \oplus \C x = \p_{\pi'}^{+}(\var{u})$.

We now show that $x=x(\underline{\lambda})$ is semisimple.
To start with, we prove that $x$ is semisimple if and only if $(\tau_2  \tau_5)/ \tau_1 + (\tau_4 \tau_6)/ \tau_3  \not= 0$.
As $\beta_2$ and $\beta_3$ are both different from $\alpha_{i_1},\alpha_{i_2}$ and $\alpha_{i_1}+\alpha_{i_2}$,
the component of $x$ on $\mathfrak{l}_{\pi'}$ in the decomposition
$\p_{\pi'}^{+}=\l_{\pi'}\oplus \mathfrak{m}_{\pi'}^{+}$ is $x_{-\eps{\pi'}} + \mu_1 x_{\eps{j_1}} + \nu x_{-\eps{j_1}}$.
By what foregoes, $\mu_1=(b  c_1)/a_{j_1} \not=0$.
Therefore, $x$ is semisimple if and only if $\nu\not=0$.
We have $\nu b \tau_0+ \lambda_2
a_{j_3}  \tau_5 + \lambda_3   a_{j_2} \tau_6  =0$, that is, by the choices of $\lambda_2$ and
$\lambda_3$:
$ \nu  b  \tau_0 - (a_{j_2}  \tau_2 a_{j_3} \tau_5)/ (a_{j_0}  \tau_1  )
- (a_{j_3}  \tau_4 a_{j_2}  \tau_6 ) / (a_{j_0}  \tau_3) =0$
Hence $\nu =  1 /(b  \tau_0)  \times (a_{j_2}  a_{j_3} )/ a_{j_0} \times
\left((\tau_2  \tau_5)/\tau_1 + (\tau_4  \tau_6 ) / \tau_3 \right)$.
As a result, $\nu \not= 0$ if and only $ (\tau_2  \tau_5)/\tau_1 + (\tau_4  \tau_6 ) / \tau_3\not=0$.

It remains to check that the condition $ (\tau_2  \tau_5)/\tau_1 + (\tau_4  \tau_6 ) / \tau_3\not=0$ holds.
We check the condition for the all cases considered in the proof of Lemma~\ref{l-2}.
Note that the computations of the integers $\tau_i$ can be done using \texttt{GAP}.

Type F$_4$: One checks that
$\tau_1=\tau_3=1$ and $\tau_2=\tau_4=\tau_5=\tau_6=-1$.

Type E$_6$: One checks that $\tau_1=\tau_2=\tau_3=\tau_4=1$ and $\tau_5=\tau_6=-1$.

Type E$_7$: One checks that $\tau_1=\tau_2=\tau_3=\tau_4=\tau_5=\tau_6=-1$.

Type E$_8$: One checks that $\tau_1=\tau_2=\tau_3=\tau_4=\tau_5=\tau_6=1$.
\end{proof}

To summarize, this gives us:

\begin{theorem} \label{t-2}
For simple $\g$ of exceptional type, and simple $\pi' \subset\pi$ of rank 2 containing $\alpha_{\pi}$,
the parabolic subalgebra $\mathfrak{p}_{\pi'}^{+}$ is quasi-reductive.
\end{theorem}

Using Theorem~\ref{t-2}, we obtain new cases
of quasi-reductive parabolic subalgebras in E$_6$:

\begin{theorem} \label{t-2E6}
For simple $\g$ of type E$_6$ and $\pi''=\{\alpha_1,\alpha_2,\alpha_4\}$
or $\{\alpha_1,\alpha_2,\alpha_4,\alpha_6\}$, $\p_{\pi''}^+$ is quasi-reductive.
\end{theorem}
Note that Theorem~\ref{t-2E6} cannot be deduced from Theorem~\ref{t-add}
even though $\pi''$ is not connected.
Indeed Theorem~\ref{t-add} fails in type E$_6$ as explained in
Remark~\ref{r-add}.

\begin{proof}
We approach the two cases in the same way.

Let $\pi'$ be the subset $\{\alpha_2,\alpha_4\}$.
Then $\pi'$ is a connected component of $\pi''$.
Hence, one can choose $u''=u(\underline{a},\underline{b})$ such that
both $\varp{u''}{\p_{\pi''}^+}$ and $\varp{u'}{\p_{\pi'}^+}$ are regular
(for $\p_{\pi''}^+$ and $\p_{\pi'}^+$ respectively)
where $u'=u(\underline{a},b_{\pi'})$.
Let $\underline{\lambda}=(\lambda_2,\lambda_3,\mu_{0},\mu_{1},\mu_{2},\mu_{3},\nu)$
be an element of $\C^7$ such that $x=x(\underline{\lambda})$ stabilizes $\varp{u'}{\p_{\pi'}^+}$ (cf.~Lemma~\ref{l2-2}).
One can readily check that $x$ belongs to $\p_{\pi''}^+(\var{u''})$, too.
On the other hand, in both cases,
the orthogonal of $E_{\pi'',\pi}$ in $\mathfrak{h}$
has dimension 1, is contained in $\p_{\pi''}^+(\var{u''})$,
and does not contain $x$.
Hence, as $x$ is semisimple (by Lemma~\ref{l2-2}),
we have found a torus a dimension 2
which is contained in $\p_{\pi''}^+(\var{u''})$.

We distinguish now the two cases:\\
Case $\pi''=\{\alpha_1,\alpha_2,\alpha_4\}$: by~(\ref{eq-ind}), $\ind \p_{\pi''}^+=2$.
Then, the above discussion shows that $\varp{u''}{\p_{\pi''}^+}$ is of reductive type.\\
Case $\pi''=\{\alpha_1,\alpha_2,\alpha_4,\alpha_6\}$: by~(\ref{eq-ind}), $\ind \p_{\pi'''}^+ =3$.
So, it suffices to provide a nonzero semisimple element in $\p_{\pi''}^+(\var{u''})$ which does not lie in the preceding torus.
We claim that the (semisimple) element
$y=(a_{K_3}/b_{\{\alpha_6\}}) x_{\alpha_1}+ (a_{K_2}/b_{\{\alpha_1\}}) x_{-\alpha_1} + (a_{K_3}/b_{\{\alpha_1\}}) x_{\alpha_6}+
(a_{K_2}/b_{\{\alpha_6\}}) x_{-\alpha_6} +x^+$ does the job, where $x^+$ is an element of $\m_{\pi''}^+$ and
where $K_i$ (for $1\leq i \leq 4$) corresponds to the highest root $\eps{i}$.
\end{proof}

\subsection{}\label{Ss-GAP}
Using the results of Sections~\ref{S-2},~\ref{S-3} and~\ref{S-4}, we are able to deal with a large number
of parabolic subalgebras.
Unfortunately, the results obtained so far do not cover all parabolic subalgebras.
There remains a small number of cases.
We consider these here.
This will complete the proof of Theorems~\ref{t-EF} and~\ref{t-E6}.

We first consider examples which do not need of the computer programme \texttt{GAP}.

It is well known that minimal parabolic subalgebras of a real simple (finite dimensional) Lie algebra
are quasi-reductive, see e.g.~\cite{Moo}.
Moreover, the complexified subalgebras give rise to quasi-reductive subalgebras of the corresponding complex simple Lie algebra.
In type F$_4$ and type E$_6$ the so-obtained parabolic subalgebras of $\g$ correspond to the subsets
$\pi'=\{\alpha_1,\alpha_2,\alpha_3\}$ and $\pi'=\{\alpha_2,\alpha_3,\alpha_4,\alpha_5\}$ of $\pi$ respectively.
As a result, we have:

\begin{proposition}\label{p-min}
{\rm (i)} If $\g$ is of type $\mathrm{F}_4$ and if $\pi'$ is  $\{\alpha_1,\alpha_2,\alpha_3\}$
then $\p_{\pi'}^{+}$ is quasi-reductive.

{\rm (ii)} If $\g$ is of type $\mathrm{E}_6$ and if $\pi'$ is $\{\alpha_2,\alpha_3,\alpha_4,\alpha_5\}$
then $\p_{\pi'}^{+}$ is quasi-reductive.
\end{proposition}

We consider now the remaining cases.
For all these cases, we are able to find
$(\underline{a},\underline{b})\in \C^{*(\kg{\pi}+\kg{\pi'})}$
such that  $\var{u (\underline{a},\underline{b})}$ is of reductive type for $\p_{\pi'}^{+}$.
We have used the computer programme \texttt{GAP} to check that the stabilizer of such a form is a torus of $\g$.
The commands we have used are presented in Appendix \ref{app1}.

\begin{proposition} \label{p-GAP}
{\rm (i)} If $\g$ is of type $\mathrm{E}_6$ and if $\pi'$ is $\{\alpha_1,\alpha_2,\alpha_3,\alpha_4,\alpha_5\}$
or $\{\alpha_2,\alpha_3,\alpha_4,\alpha_5,\alpha_6\}$ then $\p_{\pi'}^+$ is quasi-reductive.

{\rm (ii)} If $\g$ is of type $\mathrm{E}_7$ and if $\pi'$ is one the subsets
$\{\alpha_1,\alpha_2,\alpha_3,\alpha_4\}$, $\{\alpha_1,\alpha_2,\alpha_3,\alpha_4,\alpha_5\}$,
$\{\alpha_1,\alpha_2,\alpha_3,\alpha_4,\alpha_5,\alpha_6\}$, $\{\alpha_1,\alpha_3,\alpha_4,\alpha_5\}$ or
$\{\alpha_1,\alpha_3,\alpha_4,\alpha_5,\alpha_6,\alpha_7\}$
then $\p_{\pi'}^+$ is quasi-reductive.

{\rm (iii)} If $\g$ is of type $\mathrm{E}_8$ and if $\pi'$ is one the subsets
$\{\alpha_5,\alpha_6,\alpha_7,\alpha_8\}$, $\{\alpha_3,\alpha_4,\alpha_5,\alpha_6,\alpha_7,\alpha_8\}$,
$\{\alpha_2,\alpha_4,\alpha_5,\alpha_6,\alpha_7,\alpha_8\}$ or
$\{\alpha_2,\alpha_3,\alpha_4,\alpha_5,\alpha_6,\alpha_7,\alpha_8\}$
then $\p_{\pi'}^+$ is quasi-reductive.
\end{proposition}

This proposition completes the proof of Theorems~\ref{t-EF} and~\ref{t-E6}; the other cases are dealt
with either in Remark~\ref{r-0}, or in Example~\ref{ex2-E6}, or in Theorems~\ref{t-one},~\ref{t-2} and \ref{t-2E6}
(or deduced from Theorem~\ref{t-tran} or Theorem~\ref{t-add} as explained before).

\begin{remark} \label{r-yak}
As noticed in the introduction, Proposition~\ref{p-GAP} can be proved without the help of \texttt{GAP}; this is done
in a joint work of the second author and O.~Yakimova,~\cite{MY} where the authors consider the {\em maximal reductive stabilizers}
of quasi-reductive parabolic subalgebras of simple Lie algebras.
\end{remark}

\appendix

\section{} \label{app1}
In this appendix, we explain how to use \texttt{GAP} to verify that for
suitable  $u=u(\underline{a},\underline{b})$
and  $\pi'$ as described in Proposition \ref{p-GAP} the linear form
$\varp{u}{\p_{\pi'}^{+}}$
is of reductive type.
We do this for the example $\g=$E$_7$ and $\pi'=\{\alpha_1,\alpha_2,\alpha_3,\alpha_4,\alpha_5\}$,
the other cases work similarly.
First, we define the simple Lie algebra \texttt{L} (= $\g$), a root system \texttt{R}
and a Chevalley Basis (\texttt{h},\texttt{x},\texttt{y}) of \texttt{L},
and then the parabolic subalgebra \texttt{P} (= $\p_{\pi'}^+$) generated by \texttt{gP};
its dimension is \texttt{dP}:

{\normalsize

\texttt{>L:=SimpleLieAlgebra("E",7,Rationals);;R:=RootSystem(L);;}

\texttt{>x:=PositiveRootVectors(R);;y:=NegativeRootVectors(R);;}

\texttt{>g:=CanonicalGenerators(R);;h:=g[3];;}

\texttt{>gP:=Concatenation(g[1],h,y$\texttt{\{[1..5]\}}$);;P:=Subalgebra(L,gP);;dP:=Dimension(P);}\\
\texttt{90}}

\noindent
Next we choose
numbers $(\texttt{a1,a2,a3,a4,a5,a6,a7,b1,b2,b3,b4})\in (\C^{*})^{(\kg{\pi}+\kg{\pi'})}$
and we define the element \texttt{u=u1+u2} (= u(\underline{a},\underline{b}))
of $\p_{\pi'}^-$:
{\normalsize

\texttt{>a1:=-3;;a2:=5;;a3:=7;;a4:=11;;a5:=13;;a6:=-17;;a7:=19;;}

\texttt{>b1:=23;;b2:=-29;;b3:=31;;b4:=37;;}

\texttt{>u2:=a1*y[63]+a2*y[49]+a3*y[28]+a4*y[7]+a5*y[2]+a6*y[3]+a7*y[5];;}

\texttt{>u1:=b1*x[37]+b2*[16]+b3*[4]+b4*x[1];;u:=u1+u2;;}

\noindent
We are now ready to compute the  stabilizer of $\varp{\texttt{u}}{\texttt{P}}$.
To start with, we calculate the vector space \texttt{V} generated by the brackets
\texttt{u*bP[i]}, for $i=1,\ldots,\texttt{dP}$, where \texttt{bP} is a basis of \texttt{P}.
We obtain the orthogonal \texttt{K} of \texttt{V} with respect to the Killing
form thanks to the command \texttt{KappaPerp}.
Then, the stabilizer \texttt{S} of $\varp{\texttt{u}}{\texttt{P}}$ is the intersection of
\texttt{K} and \texttt{P}:

{\normalsize

\texttt{>bP:=List(Basis(P));; l:=[];;for i in [1..dP] do l[i]:=u*bP[i];od;;l;}

\texttt{>V:=Subspace(L,l);;K:=KappaPerp(L,V);;S:=Intersection(K,P);;dS:=Dimension(S);}\\
\texttt\texttt{4}}

\noindent
The fact $\dim\,$\texttt{S}=4 shows that $\varp{\texttt{u}}{\texttt{P}}$ is regular,
since $\ind \texttt{P}=4$.
It remains to check that \texttt{S} is a reductive subalgebra of \texttt{L}..
To process, we check that the restriction of the Killing form to $\texttt{S}\times \texttt{S}$ is nondegenerate. For that
it suffices to compute the intersection between \texttt{S} and its orthogonal
in \texttt{L}.
The result has to be a vector space of dimension 0:

{\normalsize

\texttt{>KS:=Intersection(KappaPerp(L,S),S);}\\
\texttt{<vector space of dimension 0 over Rationals>}}



\begin{thebibliography}{DK00}

\bibitem[B02]{Bou} N.~Bourbaki, {\em Lie groups and Lie algebras. Chapters 4--6}, Elements of Mathematics (Berlin), Springer-Verlag, Berlin, 2002.

\bibitem[DK00]{DK} V.~Dergachev, A.A.~Kirillov, {\it Index of Lie algebras of seaweed type}, J. Lie Theory,  {\bf 10},  331--343 (2000).

\bibitem[Du82]{Du1} M.~Duflo, {\it Th\'eorie de Mackey pour les groupes de Lie alg\'ebriques}, Acta Math. 149 (1982), no. 3-4, 153--213.

\bibitem[DKT]{DKT} M.~Duflo, M.S.~Khalgui and P.~Torasso, {\it Quasi-reductive Lie algebras}, preprint.

\bibitem[DV69]{DuV} M.~Duflo and M.~Vergne, {\it Une propri\'et\'e de la repr\'esentation coadjointe d'une alg\`ebre de {L}ie}, C.~R.~Acad.~Sci. Paris S\'er. A-B, {\bf 268}, A583--A585, (1969).

\bibitem[Dv03]{Dv} A.~Dvorsky, {\it Index of parabolic and seaweed subalgebras of {$\mathfrak{so}_n$}},  Lin. Alg. Appl, {\bf 374}, 127--142 (2003).

\bibitem[El]{El} A.G.~Elashvili, {\it On the index of parabolic subalgebras of semisimple Lie algebras}, communicated by M. Duflo.

\bibitem[HC65]{HC1} Harish-Chandra, {\it Discrete series for semisimple Lie groups, I. Construction of invariant eigendistributions}, Acta Math. {\bf 113}, 241--318 (1965).

\bibitem[HC66]{HC2} Harish-Chandra, {\it Discrete series for semisimple Lie groups, II. Explicit determination of the characters}, Acta Math. {\bf 116}, 1--111 (1966).

\bibitem[J06]{Jo1} A.~Joseph, {\it On semi-invariants and index for biparabolic (seaweed) algebras. {I}}, J. Algebra,  {\bf 305}, no.~1,
485--515 (2006).

\bibitem[J07]{Jo2} A.~Joseph, {\it On semi-invariants and index for biparabolic (seaweed) algebras. {II}}, J. Algebra,  {\bf 312}, no.~1,
158--193 (2007).

\bibitem[Ki68]{Kir} A.A.~Kirillov, {\em The method of orbits in the theory of unitary representations of Lie groups}, Funkcional. Anal. i Prilozen. no.~1, 96--98 (1968).

\bibitem[Moo70]{Moo} C.C.~Moore, {\it Restrictions of unitary representations to subgroupes and ergodic theory:
Group extensions and group cohomology}, Lectures Notes in Physics, {\bf 6} (1970).

\bibitem[Mor06a]{Mo1} A.~Moreau, {\it Indice du normalisateur du centralisateur d'un \'el\'ement nilpotent dans une alg\`ebre de {L}ie semi-simple},
Bull. Soc. Math. France {\bf 134}, no.~1, 83--117 (2006).

\bibitem[Mor06b]{Mo2} A.~Moreau, {\it Indice et d\'ecomposition de Cartan d'une alg\`ebre de Lie semi-simple r\'eelle}, J. Algebra, {\bf 303},
 no.~1, 382--406 (2006).

\bibitem[MY]{MY} A.~Moreau and O.~Yakimova, {\it Coadjoint orbits of reductive type of seaweed algebras},
preprint (2010).

\bibitem[P01]{Pa1} D.I.~Panyushev, {\it Inductive formulas for the index of seaweed Lie algebras},  Moscow Math. J., \textbf{2}, {\bf 1}, 221--241, 303 (2001).

\bibitem[P03]{Pa2} D.I.~Panyushev, {The index of a Lie algebra, the centraliser of a nilpotent element, and the normaliser of the centraliser},
Math. Proc. Cambridge Philos. Soc. {\bf 134}, no.~1,  41--59 (2003).

\bibitem[P05]{Pa3} D.I.~Panyushev, {\it An extension of Ra\"{i}s' theorem and seaweed subalgebras of simple Lie algebras},
Ann. Inst. Fourier (Grenoble), {\bf 55}, no.~3, 693--715 (2005).

\bibitem[TY04a]{ty0} P.~Tauvel and R.W.T.~Yu, {\it Indice et formes lin\'eaires stables dans les alg\`ebres de Lie}, J. Algebra, {\bf 273}, 507--516 (2004).

\bibitem[TY04b]{TY2} P.~Tauvel and R.W.T.~Yu, {\it Sur l'indice de certaines alg\`ebres de Lie}, Ann.Inst. Fourier (Grenoble) {\bf 54}, 1793--1810 (2005).

\bibitem[TY05]{TY1} P.~Tauvel and R. W.T.~Yu, {\it Lie Algebras and Algebraic groups}, Monographs in Mathematics (2005), Springer, Berlin Heidelberg New York.

\bibitem[Ti62]{Ti1} J.~Tits, {\it Groupes semi-simple isotropes}, Colloque sur la th\'eorie des groupes alg\'ebriques, C.B.R.M., Bruxelles, 137--147 (1962).

\bibitem[Ti66]{Ti2} J.~Tits, {\it Classification of algebraic semisimple groups}, Proc. Symposia Pure Math. AMS, 33--62 (1966).

\bibitem[Ya06]{Ya} O.~Yakimova, {\it The index of centralizers of elements in classical Lie algebras}, Funktsional. Anal. i Prilozhen. {\bf 40}, no.~1, 52--64, 96 (2006).

\end{thebibliography}
\end{document}